\documentclass[12pt,a4paper]{amsart}
\usepackage[english]{babel}
\theoremstyle{plain}
\usepackage{amssymb,amsfonts,mathrsfs,amscd}

\advance\hoffset-20mm \advance\textwidth40mm

\numberwithin{equation}{section}

\numberwithin{equation}{section}
\theoremstyle{plain}
\newtheorem{theorem}{Theorem}[section]
\newtheorem{lemma}{Lemma}[section]

\newtheorem{proposition}{Proposition}[section]

\theoremstyle{definition}
\newtheorem{definition}{Definition}[section]
\newtheorem{remark}{Remark}[section]
\newtheorem{example}{Example}[section]

\theoremstyle{definition}

\theoremstyle{definition}

\def\AA{{\mathbb A}}
\def\XX{{\mathbb X}}
\def\NN{{\mathbb N}}
\def\ZZ{{\mathbb Z}}
\def\QQ{{\mathbb Q}}
\def\TT{{\mathbb T}}
\def\kk{{\Bbbk}}
\def\CC{{\mathbb C}}
\def\RR{{\mathbb R}}
\def\PP{{\mathbb P}}
\def\Of{{\mathscr{O}}}
\def\Lf{{\mathscr{T}}}
\def\fg{{\mathfrak g}}
\def\fh{{\mathfrak h}}

\def\cone{\mathop{\rm Cone}}
\def\span{\mathop{\rm span}}
\def\Aut{\mathop{\rm Aut}}
\def\SAut{\mathop{\rm SAut}}

\def\Stab{\mathop{\rm Stab}}
\def\Spec{\mathop{\rm Spec}}
\def\AC{\mathop{\rm AffCone}}
\def\Der{\mathop{\rm Der}}

\def\deg{\mathop{\rm deg}}

\def\reg{\mathop{\rm reg}}
\def\ord{\mathop{\rm ord}}
\def\rk{\mathop{\rm rk}}
\def\Char{\mathop{\rm Char}}

\def\ML{\mathop{\rm ML}}
\def\FML{\mathop{\rm FML}}
\def\Pic{\mathop{\rm Pic}}
\def\Int{\mathop{\rm Int}}
\def\SP{\mathop{\rm Susp}}
\def\card{\mathop{\rm card}}

\def\ll1{l_{\lambda}^{-1}(1)}
\def\lm1{l_{\mu}^{-1}(1)}
\newcommand\qmatrix[2][1]{\left(\renewcommand\arraystretch{#1}
\begin{array}{*{20}r}#2\end{array}\right)}

\begin{document}
\sloppy

\title[Three instances of infinite transitivity]{Flag varieties,
toric varieties, and suspensions:\\Three instances of infinite
transitivity}

\author{I.\,V.~Arzhantsev}
\address{Department of Algebra, Faculty of Mechanics and Mathematics,
Moscow State University, Leninskie Gory 1, GSP-1, Moscow, 119991,
Russia}
\email{arjantse@mccme.ru}
\author{K.\,G.~Kuyumzhiyan}
\address{National Research University Higher School of Economics,\\
Laboratory of Algebraic Geometry and Its Applications,\\
7 Vavilova Str., Moscow, 117312, Russia}
\email{karina@mccme.ru}
\author{M.\,G.~Zaidenberg}
\address{Universit\'e Grenoble I, Institut Fourier, UMR 5582
CNRS-UJF, BP 74, 38402 St. Martin d'H\`eres c\'edex, France}
\email{Mikhail.Zaidenberg@ujf-grenoble.fr}

%


\begin{abstract}
We say that a group~$G$ acts
infinitely transitively  on a set~$X$ if for every $m\in\NN$ the
induced diagonal action of~$G$ is transitive on the cartesian $m$th
power $X^m\setminus\Delta$ with the diagonals removed.
We describe three classes of
affine algebraic varieties such that their automorphism groups
act infinitely transitively on their smooth loci.
The first class consists of normal affine cones over flag varieties, the second of
non-degenerate affine toric varieties, and the third of iterated suspensions over
affine varieties with infinitely transitive automorphism groups of
a reinforced type.
%
\end{abstract}




\thanks{The first author was supported by the Deligne's 2004 Balzan prize
in Mathematics. The work of the first two authors was partially
supported by Dmitry Zimin fund ``Dynasty'', the second author is partially supported by
AG Laboratory NRU-HSE, RF government
grant, ag. 11.G34.31.0023, ``EADS Foundation Chair in Mathematics'' and
Russian-French Poncelet Laboratory (UMI 2615
CNRS). The main part of the work was done
during a stay of the first author at the Institut Fourier,
Grenoble. The authors thank Institute Fourier for its generous
support and hospitality.}

\maketitle

\section*{Introduction}

All varieties in this paper are
assumed being reduced and irreducible. Unless we explicitly
precise a setting, the ground field $\kk$ is supposed to be
algebraically closed of characteristic zero.

An effective action of the additive group ${\mathbb G}_{\rm a}(\kk)$ on an
algebraic variety~$X$ defines a one-parameter unipotent subgroup
of the automorphism group $\Aut(X)$. We let $\SAut(X)$ denote the
subgroup of $\Aut(X)$ generated by all its one-parameter unipotent
subgroups. In the sequel, we adopt the following definitions.

\begin{definition}\label{def1}
Let $X$ be an algebraic variety over $\kk$. We say that a point
$x\in X$ is \textit{flexible} if the tangent space $T_x X$ is
spanned by the tangent vectors to the orbits $H.x$ of
one-parameter unipotent subgroups $H\subseteq\Aut (X)$. The
variety~$X$ is called \textit{flexible} if every smooth point
$x\in X_{\rm reg}$ is.
\end{definition}

\begin{definition}\label{def2}
An action of a group $G$ on a set $A$ is said to be
\textit{$m$-transitive} if for every two tuples of pairwise
distinct points $(a_1,a_2,\ldots,a_m)$ and $(b_1,b_2,\ldots,b_m)$
in~$A$ there exists $g\in G$ such that $g(a_i)=b_i$,
$i=1,2,\dots,m$. The action which is $m$-transitive for all
$m\in\NN$ will be called \textit{infinitely transitive}.
\end{definition}

Clearly, for $m>\dim G$ the Lie group~$G$ cannot act
$m$-transitively on a variety $X$. In fact, due to
the following theorem it cannot
act even $4$-transitively on a smooth simply connected variety.

\begin{theorem}\label{borel} ({\rm A.\ Borel}
\cite[Th\'eor\`emes 5-6]{Bo})
There exist no $3$-transitive actions of
a real Lie group $G$ on a simply connected non-compact variety,
and such $2$-transitive actions exist only on the Euclidean spaces
$\RR^n$ for $n\ge 2$. There exist no $4$-transitive actions of $G$ on a
compact simply connected variety, and such $3$-transitive actions
exist only on the spheres $S^n$ for $n\ge 2$.
\end{theorem}

For a classification of $2$- and
$3$-transitive Lie group actions see, e.g., \cite{Kr} and~\cite{Ti}.
See also~\cite{Po3} for the "generic transitivity" of an action of an algebraic group $G$,
i. e. the existence of Zariski open orbits in the Cartesian powers of a given algebraic $G$-variety.

The infinite transitivity of the full automorphism group is well
known for the affine space $\AA^n$ over~$\kk$, where $n\ge 2$; see
also~\cite{AL} and~\cite{RR} for the analytic counterpart and generalizations.
This
phenomenon takes place as well for any hypersurface $X$ in
$\AA^{n+1}$ given by the equation $uv-f(x_1,\ldots,x_{n-1})=0$, where
$n\ge 3$ and $f\in\kk[x_1,\ldots,x_{n-1}]$ is a non-constant
polynomial \cite[\S 5]{KZ}. In the sequel we call  such a variety~$X$
a \textit{suspension} over $Y=\AA^{n-1}$. In the analytic
category, similar phenomenons were studied in the spirit of the
Andersen-Lempert-Varolin theory~\cite{AL}, \cite{Va}, see, e.g.,
\cite{Fo}, \cite{KK1}, \cite{Ro}, and \cite{TV}. This concerns, in
particular, the infinite transitivity on smooth affine algebraic
varieties of certain subgroups of biholomorphic transformations
generated by complete regular vector fields, see survey~\cite{KK2},
especially \S 2(B) and Remark~2.2.

Following~\cite[\S 5]{KZ}, we are interested here in algebraic
varieties~$X$ such that the special automorphism group $\SAut (X)$
acts infinitely transitively on the smooth locus $X_{\rm reg}$ of~$X$.
We show, in particular, that this property is preserved when
passing to a suspension.

\begin{definition}\label{def3}
By a \textit{suspension} over an affine variety~$Y$ we mean a
hypersurface $X\subseteq Y\times \AA^2$ given by the equation
$uv-f(y)=0$, where $\AA^2=\Spec\kk[u,v]$ and $f\in\kk[Y]$ is
non-constant. In particular, $\dim X=1+\dim Y$.
\end{definition}

Our main results can be formulated as
follows. An affine variety is called {\it non-degenerate} if every invertible regular
function on it is constant.

\begin{theorem}\label{main0}
\begin{enumerate}
\item Consider a flag variety $G/P$, where~$G$ is a semisimple
algebraic group, and~$P$ is a parabolic subgroup in~$G$.
Then every normal affine cone~$X$ over $G/P$ is flexible and its
special automorphism group $\SAut (X)$ acts infinitely
transitively on the smooth locus $X_{\rm reg}$.
\item The analogous
conclusion holds if~$X$ is any non-degenerate
affine toric
variety 
of dimension at least~2.
\item Suppose that an
affine variety~$X$ 
is flexible and either $X=\AA^1$,
or $\dim X\ge 2$ and the special automorphism group $\SAut (X)$
acts infinitely transitively on the smooth locus $X_{\rm reg}$.
Then all iterated suspensions over~$X$ also have
properties of flexibility and infinite transitivity of
the special automorphism group.
\end{enumerate}
\end{theorem}

Theorem~\ref{main0}(n) is proved in Section~$n$, where $n=1,2,3$,
respectively.

Smooth compact real algebraic surfaces with
infinitely transitive automorphism groups were classified
in~\cite{BH}, \cite{BM}, \cite{HM1}, and~\cite{HM2}. In
Theorem~\ref{main10} we extend part~3 of Theorem~\ref{main0} to
real algebraic varieties, under certain additional restrictions.
These restrictions were weakened in a recent paper~\cite{KM}.
In the particular case of a suspension over the affine line, the result
remains valid over an arbitrary field of characteristic zero, see
Theorem~\ref{aa}.

Recall~\cite[\S 9]{Fre} that the Makar-Limanov invariant $\ML(X)$
of an affine variety~$X$ is the intersection of the kernels of all
locally nilpotent derivations of $\kk[X]$, or, in other words, the
subalgebra in $\kk[X]$ of common invariants for all one-parameter
unipotent subgroups of $\Aut (X)$. From this definition it is
straightforward that $\ML(X)=\kk[X]^{\SAut(X)}$. Hence the
Makar-Limanov invariant of~$X$ is trivial (that is $\ML(X)=\kk$)
provided that the special automorphism group $\SAut(X)$ acts on~$X$
with a dense open orbit (cf.~\cite{Po2}). In particular, this holds
for the varieties in all the three classes from Theorem~\ref{main0}
(for the first two of them, see also~\cite[3.16]{KPZ}, \cite{Li}, and~\cite{Po2}). On the other hand, $\ML(X)$ is trivial
if~$X$ is flexible. Indeed, if $f\in\kk[X]^{\SAut(X)}$ then the
differential~$df$ vanishes along the orbits of any unipotent
subgroup, hence it vanishes on the tangent space at any flexible
point of $X_{\reg}$. Since~$X$ is flexible, $f$ is constant.

In the first preprint version of this paper we conjectured
that for a variety of dimension at least 2 over an algebraically
closed field of characteristic zero {\em the transitivity  of the
group $\SAut(X)$ on $X_{\rm reg}$, the infinite transitivity of
this group, and flexibility of~$X$ are equivalent}. Later on this
conjecture was proved in~\cite{AFKKZ}.
Furthermore, in the same preprint version we proposed to {\em
characterize the class of affine varieties such that the group
$\SAut(X)$ acts infinitely transitively on a Zariski open subset.}
This was also done in~\cite{AFKKZ} (see Theorem~2.2 and
Proposition~5.1) in terms of the field Makar-Limanov invariant
$\FML(X)$, introduced by A.~Liendo~\cite{Lie2}. In~\cite{Perep},
it is shown that affine cones over del Pezzo surfaces of degrees~4
and~5 are flexible, which implies infinite transitivity of the
action of $\SAut(X)$ on these varieties.

We are grateful to Dmitry Akhiezer, Shulim Kaliman, and Alvaro
Liendo for useful discussions and references.

\section{Affine cones over flag varieties}
Given a
connected, simply connected, semisimple algebraic group $G$
over~$\kk$, we consider an irreducible representation $V(\lambda)$
of~$G$ with a highest weight~$\lambda$ and a highest weight vector
$v_0\in V(\lambda)$. Let
$$
Y=G[v_0]\subseteq \PP(V(\lambda))
$$
be the closed $G$-orbit
of $[v_0]$ in the associate projective representation, and
$$
X=\AC(Y)=\overline{Gv_0}=Gv_0\cup\{0\}
$$
be the affine cone over~$Y$. Note that such a cone~$X$ is called
an {\it HV-variety} in terminology of~\cite{PV}.

\medskip

\begin{remark}\label{prnor}
Actually every projective embedding
$\varphi:G/P\stackrel{\cong}{\longrightarrow} Y\hookrightarrow\PP^n$
with a projectively normal image~$Y$, where $P\subseteq G$ is a
parabolic subgroup, arises in this way. Indeed, being projectively
normal, $Y$ is as well linearly normal, i.e.,
$\varphi=\varphi_{|D|}$, where $D\in\Pic (G/P)$ is very ample.
Hence $D\thicksim \sum_{i=1}^s a_iD_i$, where $D_1,\ldots,D_s$ are
the Schubert divisorial cycles on $G/P$, and $a_i\in\ZZ$, $a_i>0$ for all
$i=1,\ldots,s$. Then $V(\lambda)=H^0(G/P, {\mathscr
O}_{G/P}(D))^*$ is a simple $G$-module with a highest weight
$\lambda=\sum_{i=1}^s a_i\omega_i$, where $\omega_1,\ldots,
\omega_s$ are fundamental weights,  and $Y=G[v_0]$ for a highest
weight vector $v_0$; see, e.g., \cite[Theorem 5]{Po1}.
\end{remark}

\begin{remark}
Recall that an \textit{end} in the Freudenthal sense of a topological variety~$M$
is the equivalence class of decreasing sequences of connected open subsets in~$M$, such that their boundaries are compact, and the intersection of complements of these sets is empty.
Let $X$ be the affine cone introduced above, and let $\kk=\CC$. Then
the homogeneous variety
$X\setminus\{0\}$ has two ends. Let $G$
be a connected Lie group and $H\subseteq G$ a closed connected
subgroup. According to A.~Borel~\cite[Th\'eor\`eme 2]{Bo}, the
homogeneous space $G/H$ has at most two ends.
D.~Akhiezer~\cite{Ak}\footnote{See also~\cite{HO} and~\cite{Le} for the
complex analytic counterpart in the context of Lie groups.}
showed that if $G$ is a linear algebraic group over ${\CC}$
and
$H\subseteq G$ is an algebraic subgroup (not necessarily connected)
then $G/H$ has exactly two
ends if and only if $H$ is the kernel of a nontrivial character
$\chi\colon P\to {\mathbb G}_{\rm m}(\CC)$, where $P$ is a parabolic subgroup of~$G$,
and ${\mathbb G}_{\rm m}(\CC)$ is the multiplicative group of the field~$\CC$.
The homogeneous fibration $G/H\rightarrow G/P$ realizes
$G/H$ as a principal ${\mathbb G}_{\rm m}(\CC)$-bundle over the homogeneous
projective variety\footnote{In this general setting,
neither the pair $(G,P)$ nor $(G,H)$ is uniquely defined by the
flag variety $G/P$. } $G/P$. Furthermore, $G/H$ admits a projective
completion by two disjoint divisors $E_0$ and $E_\infty$, where
$E_0\cong E_\infty\cong G/P$, and $\tilde X:=G/H\cup E_0\to G/P$
represents a line bundle, say, $L$ over $G/P$. Its zero section
$E_0$ is contractible if and only if the dual line bundle $L^{-1}$
is ample, and then also very ample. In the latter case the
contraction of $E_0$ yields the affine cone $X$ over the image
$Y=\varphi_{|L^{-1}|}(G/P)\subseteq\PP^n$. For $\kk=\CC$, any
affine cone~$X$ as in Theorem~\ref{flag} arises in this way, with
$H\subseteq P$ being the stabilizer of a smooth point of the cone~$X$,
$\tilde X$ the blowup of~$X$ at the vertex, and~$E_0$ the
exceptional divisor.
\end{remark}

The next result provides part~1 of Theorem~\ref{main0}.

\begin{theorem}\label{flag}
Let  $X$ be the affine cone
over a flag variety $G/P$ under an embedding $G/P\hookrightarrow
\PP^N$ with projectively normal image. Then~$X$ is flexible and
the group $\SAut(X)$ acts $m$-transitively on $X\backslash \{0\}$
for any $m\in\NN$.
\end{theorem}

The flexibility follows from the next general observation.

\begin{proposition}\label{flag0}
If a semisimple linear
algebraic group~$G$ acts on an affine variety~$X$ and this action
is transitive on $X_{\rm reg}$, then~$X$ is flexible.
\end{proposition}

\begin{proof}
The group~$G$ acts on~$X$
with an open orbit $X_{\rm reg}=G.x_0$.
The dominant morphism  onto this orbit
$\varphi:G\to X$, $g\longmapsto g.x_0$, yields a surjection
$d\varphi:\fg\to T_{x_0}X$, where $\fg={\rm Lie}(G)$. We claim
that $\fg$ is spanned by nilpotent elements over $\kk$, which
implies the assertion. Indeed, consider the decomposition
$\fg=\bigoplus_{i=1}^k\fg_i$ of $\fg$ into simple ideals. 
 Let
$\fh$ be the span of the set of all nilpotent elements in $\fg$.
This is an ad-submodule of $\fg$ and so an ideal of $\fg$, hence a
direct sum of some of the simple ideals $\fg_i$. However every
simple ideal $\fg_i$, $i=1,\ldots,k$, contains at least one nonzero
nilpotent element. Therefore $\fh=\fg$, as claimed.
\end{proof}

In the setting of Theorem~\ref{flag}, $X_{\rm reg}=X$ if
$X\cong\AA^n$ and $X_{\rm reg}=X\setminus\{0\}$ otherwise. Anyhow,
the group~$G$ acts transitively on $X\setminus\{0\}$, see, e.g.,
\cite[Theorem 1]{PV}. Hence by Proposition~\ref{flag0} $X$ is
flexible.

Before passing to the proof of infinite transitivity we need some
preparation.

Let $P\subseteq G$ be the stabilizer of the line $\langle
v_0\rangle\subseteq V(\lambda)$, $B=TB_u\subseteq P$ be a Borel
subgroup of~$G$ with the maximal torus~$T$ and the unipotent
radical $B_u$, and $\XX(T)$ be the character lattice of~$T$.
Consider the weight decomposition
$$
V(\lambda)=\bigoplus_{\nu\in \XX(T)} V(\lambda)_{\nu}=\langle
v_0\rangle\oplus H({\lambda})\,,
$$
where $\langle v_0\rangle=V(\lambda)_\lambda\,$, and
$H({\lambda})\subseteq V(\lambda)$ is the hyperplane
$$
H(\lambda)=\bigoplus_{\nu\in \XX(T)\setminus\{\lambda\}}V(\lambda)_{\nu}\,.
$$
The coordinate
function $l_{\lambda}\in V(\lambda)^*$ of the first projection
$p_1:v\longmapsto l_{\lambda}(v)v_0$ defines a non-trivial
character of~$P$.

Let $B^-=TB_u^-$ be the Borel subgroup of~$G$ opposite to $B=B^+$.
The flag variety $G/P$ contains an open $B^-$-orbit (the big
Schubert cell) isomorphic to the affine space $\AA^n$, where
$n=\dim G/P$. Its complement is a union of the divisorial Schubert
cycles $D_1,D_2,\ldots, D_s$, see e.g.~\cite[pp. 22--24]{LR}.

The orbit map $G\to \PP (V(\lambda)),\, g\mapsto g. [v_0]$, embeds
$G/P$ onto a subvariety $Y\subseteq \PP (V(\lambda))$. Let
$\omega_\lambda\subseteq Y$ be the image of the big Schubert cell
under this embedding. By~\cite[Theorem~2]{Po1} the hyperplane
$$
\mathscr H(\lambda)=\PP
(H(\lambda))=l_\lambda^{-1}(0)\subseteq\PP(V(\lambda))\,
$$
is supported by the union of the Schubert divisors $\bigcup_{i=1}^s
D_i$. In particular
$\omega_\lambda=Y\setminus \mathscr
H(\lambda)$.

Let $\sigma \colon \widehat X\to X$ be the blow-up of the cone~$X$
at the vertex~$0$. The exceptional divisor $E\subseteq \widehat X$
is isomorphic to~$Y$. Moreover, the natural map $\pi\colon
X\setminus \{0\}\to Y$ yields the projection $p\colon \widehat
X\to Y$ of the line bundle $\mathscr O_Y(-1)$ on $Y$ with~$E$
being the zero section. Since $\omega_\lambda\cong \AA^n$, the
restriction of $\mathscr O_Y(-1)$ to $\omega_\lambda$ is
a trivial line bundle. Hence  the open set
$$
\Omega_\lambda:= \pi^{-1}(\omega_\lambda)=X\setminus
H(\lambda)\subseteq X\setminus
\{0\}\cong \widehat X\setminus E
$$
is isomorphic to $\AA^n\times \AA^1_*$, where $\AA^1_*=\AA^1\setminus\{0\}$.

For every $c\in\AA^1_*$ the invertible function $l_\lambda(\cdot,
c)$  is constant on the affine space~$\AA^n$. Thus on
$\Omega_\lambda$, which is isomorphic to $\AA^n\times \AA^1_*$, we have
$l_\lambda=az^k$ for some $a\in {\mathbb G}_{\rm m}(\kk)$, where $z$
is a coordinate in $\AA^1_*$. Actually, here $k=1$, since
$l_\lambda$ gives a coordinate on $\langle v_0\rangle$. We may also
assume that $a=1$ and so
$l_\lambda\vert\,_{\Omega_\lambda}:\Omega_\lambda\to \AA^1_*$ is
the second projection.

To prove the infinite transitivity  of the group $\SAut (X)$ on
$X\setminus \{0\}$ as stated in Theorem~\ref{flag}, let us first show
the infinite  transitivity of $\SAut (X)$ on each hyperplane
section $\Omega_\lambda(c_0):=l_\lambda^{-1}(c_0)\subseteq X$,
where $c_0\ne 0$; cf.~\cite[Lemma 5.6]{KZ}. More precisely, given
a $k$-tuple of distinct points $c_1,\ldots,c_k\in\kk$ different
from $c_0$, we consider the subgroup
$\Stab^\lambda_{c_1,\ldots,c_k}\subseteq\SAut(X)$ of all
automorphisms fixing pointwise the subvarieties
$\Omega_\lambda(c_i)$ for all $i=1,\ldots,k$ and leaving invariant
the function $l_\lambda$.

\begin{proposition}\label{stab}
In the notation as above,
for every $n\ge 2$ and every $m\in \NN$ the group $\Stab^\lambda_{c_1,\ldots,c_k}$ acts
$m$-transitively on $\Omega_\lambda(c_0)\cong \AA^n$.
\end{proposition}

\begin{proof}
Let $Q_1,Q_2,\ldots,Q_m$ and $Q'_1,Q'_2,\ldots,Q'_m$ be two tuples
of pairwise distinct points in $\Omega_\lambda(c_0)$. For any
$n\geqslant 2$ the group $\SAut(\AA^n)$ acts $m$-transitively  on
$\AA^n$; see e.g.~\cite[Lemma 5.5]{KZ}. Since
$\Omega_\lambda(c_0)\cong \AA^n$, we can find $g\in\SAut
(\Omega_\lambda(c_0))$ mapping $(Q_1,Q_2,\ldots,Q_m)$ to
$(Q'_1,Q'_2,\ldots,Q'_m)$. By definition,
$g=\delta_1(1)\delta_2(1)\ldots\delta_s(1)$ for some one-parameter
unipotent subgroups $\delta_1, \delta_2,\ldots,
\delta_s\subseteq\SAut(\Omega_\lambda(c_0))$. Let
$\partial_1,\partial_2,\ldots,\partial_s$ be the corresponding
locally nilpotent derivations\footnote{A derivation $\partial$ of
a ring~$A$ is called locally nilpotent if $\forall a\in A$,
$\partial^n a=0$ for some $n\in \NN$.} (LNDs for short) of the algebra
$\kk[\Omega_\lambda(c_0)]$. First we extend
them to LNDs
$\overline\partial_1,\overline\partial_2,\ldots,\overline\partial_s$
of $\kk[\Omega_\lambda]\cong \kk[\AA^n\times \AA^1_*]$ by putting
$\overline\partial_i(l_\lambda)=0$.

Recall that $\Omega_\lambda$ is a principal Zariski open subset
in~$X$ defined by the function~$l_\lambda$. In particular, for every
$i=1,\ldots,s$ we have $\overline\partial_i\colon \kk[X]\to
\kk[X][1/{l_\lambda}]$. Since $\kk[X]$ is finitely generated,
there exists $N\in\NN$ such that $(l_\lambda)^N\, \overline
\partial_i$ is an LND of $\kk[X]$ for all $i=1,\ldots,s$;
cf.~\cite[Proposition ~3.5]{KPZ}.

Let  $q[z]\in\kk[z]$ be a polynomial with $q(c_0)=1$ which has
simple roots at $c_1,\ldots,c_k$ and a root $z=0$ of multiplicity
$N$ (we recall that $c_0\neq 0$). Then for every $i=1,\ldots,s$,
$q(l_\lambda)\overline\partial_i$ is an LND of $\kk[X]$ such that
the corresponding one-parameter subgroup in
$\Stab^\lambda_{c_1,\ldots,c_k}$ extends the subgroup $\delta_i$.
Thus~$g$ extends to an element of the group
$\Stab^\lambda_{c_1,\ldots,c_k}$. Now the assertion follows.
\end{proof}

Let  $\mu$ be an extremal weight of the simple $G$-module
$V(\lambda)$ different from $\lambda$. Then~$\mu$ defines a
parabolic subgroup $P'$ conjugated to $P$, the corresponding
linear form $l_\mu\in V(\lambda)^*$, and the principal Zariski
open subset $\Omega_\mu=\{l_\mu\neq 0\}$ of~$X$, where
$X\setminus\Omega_\mu=H(\mu):=l_\mu^{-1}(0)$.

\begin{lemma}\label{not0}
For every set of $m$ distinct points
$Q_1,Q_2,\ldots, Q_m\in X\setminus \{0\}$ there exists
$g\in\SAut(X)$ such that $g(Q_i)\in \Omega_\mu$ for all
$i=1,\ldots, m$.
\end{lemma}

\begin{proof}
Since the group $G$ is semisimple, it is contained in $\SAut(X)$,
see~\cite[Lemma 1.1]{Po2}. Clearly, $G_i:=\{g\in G \,|\,
g(Q_i)\in H(\mu)\}$, $i=1,\ldots,m$, are proper closed subsets
of~$G$. Hence the conclusion of the lemma holds for every $g\in
G\setminus (G_1\cup\ldots\cup G_m)$.
\end{proof}

\begin{lemma}\label{nonconstant}
For every $c\neq 0$ the restriction
$l_\lambda\vert\,{\Omega_\mu(c)}$ is non-constant.
\end{lemma}

\begin{proof}
If the restriction $l_\lambda\vert\,{\Omega_\mu(c)}$ were a
constant equal, say, to $b$, then the cone~$X$ would be contained in
the hyperplane $bl_\mu-cl_\lambda=0$ in $V(\lambda)$, which is not
the case.
\end{proof}

{\sc Proof of Theorem~\ref{flag}.}
If $n=\dim G/P=1$ and so $G/P\cong \PP^1$, then $X$ is a normal
affine toric surface (a Veronese cone). The infinite transitivity
in this case follows from Theorem~\ref{T2} in \S~2 below.

From now on we suppose that $n\geqslant 2$. Given $m\in \NN$, we
fix an $m$-tuple of pairwise distinct points $Q_{10}, Q_{20},\ldots, Q_{m0}\in
\Omega_\lambda(1)$. Let us show that for every $m$-tuple of pairwise distinct
points $Q_1,Q_2,\ldots,Q_m\in X\setminus \{0\}$ there exists
$\psi\in\SAut(X)$ such that $\psi(Q_1)=Q_{10}, \ldots$, $\psi(Q_m)=Q_{m0}$.

According to Lemma~\ref{not0} we may suppose that $Q_i\in
\Omega_\mu$ for all $i=1,\ldots,m$. Split the set
$\{Q_1,Q_2,\ldots,Q_m\}$ into several pieces according to the
values of $l_\mu(Q_i)$:
$$
\{Q_1,Q_2,\ldots, Q_m\} \,=\,\bigsqcup _{j=1}^{k}M_j, \qquad
M_j\,=\,\{Q_i\;|\; Q_i\in\Omega_\mu(c_j)\}\,,
$$
where $c_1,\ldots,c_k\in \AA^1_*$
are pairwise distinct. By Lemma~\ref{nonconstant}, every intersection
$\Omega_\lambda(1)\cap\Omega_\mu(c_i)$ contains infinitely many
points. Acting with the subgroups $\Stab^\mu_{c_1,\ldots,\widehat
c_i,\ldots,c_k}\subseteq\SAut (X)$ (see Proposition~\ref{stab}),
we can successively send the pieces $M_i$, $i=1,\ldots,k$, to the
affine hyperplane section $\Omega_\lambda(1)$. The resulting
$m$-tuple can be sent further to the standard one
$(Q_{10},Q_{20},\ldots,Q_{m0})$ using an automorphism from
Proposition~\ref{stab} with $c_0=1$ and $k=0$. Now the proof is completed.


\section{Automorphisms of affine toric varieties}
As before, $\kk$ stands for an algebraically closed field  of
characteristic zero. In this section we consider an affine toric
variety~$X$ over~$\kk$ with a  torus $\TT$ acting effectively on~$X$.
We assume that $X$ is non-degenerate, i.e., the only
invertible regular functions on~$X$ are constants   or, which is
equivalent, that $X\not\cong Y\times\AA^1_*$, where
$\AA^1_*=\Spec\kk[t,t^{-1}]\cong {\mathbb G}_{\rm m}(\kk)$.

The following result yields part~2 of Theorem~\ref{main0}.

\begin{theorem} \label{T2}
Every non-degenerate affine toric variety~$X$ of dimension
$n\ge 1$ is flexible. If $n\ge 2$, for any $m\in\NN$ the group $\SAut(X)$ acts $m$-transitively
on the smooth locus $X_{\reg}$ of~$X$.
\end{theorem}

We note that~$X$ is flexible if $\Aut (X)$ acts transitively on
$X_{\rm reg}$ and at least one smooth point is flexible on~$X$.
Both of these properties will be established below. Let us first recall
some necessary generalities on toric varieties.

\medskip

1. {\sl Ray generators.}
Let $N$ be the lattice of one-parameter subgroups of the torus
$\TT$, $M = \XX(\TT)$ the dual lattice of characters, and
$\langle\cdot,\cdot\rangle: N\times M \to \ZZ$ the natural
pairing. Let $\chi^m$ denote the character of $\TT$ which
corresponds to a lattice point $m\in M$. We have
$\chi^m\chi^{m'}\, =\, \chi^{m+m'}$, thus the group algebra
$$
\kk[M] \ : = \ \bigoplus_{m \in M} \kk\chi^m
$$
can be identified with the algebra
$\kk[\TT]$ of regular functions on the torus $\TT$. Let $\TT . x_0$ be
the open $\TT$-orbit on~$X$. Since the orbit map $\TT \to X$, $t
\longmapsto t . x_0$, is dominant, we may identify $\kk[X]$ with a
subalgebra of $\kk[M]$. More precisely, there exists a convex
polyhedral cone $\sigma^\vee \subseteq M_{\QQ} := M\otimes_{\ZZ}
\QQ$ such that $\kk[X]$ coincides with the semigroup algebra of
$\sigma^\vee \cap M$, i.e.,
\begin{equation}
\kk[X] \ = \ \bigoplus_{m \in \sigma^\vee \cap M} \kk\chi^m\,,
\end{equation}
see~\cite{Fu} for details. We let $\sigma$
denote the cone dual to $\sigma^\vee$. The cone $\sigma$ is
pointed and of full dimension in $N_{\QQ}$. Let
$\Xi=\{\rho_1,\ldots,\rho_r\}$ be the set of ray generators,
i.e., the primitive vectors on extremal rays of the cone $\sigma$.
Given a ray generator $\rho\in \Xi$, we let $R_\rho$ denote the
associate one-parameter subgroup of $\TT$.

\medskip

2. {\sl The Orbit-Cone correspondence (see \cite[\S 3.2]{CLS}).}
There exist two natural one-to-one correspondences
$\delta\stackrel{1\div 1}{\longleftrightarrow} \tau\stackrel{1\div
1}{\longleftrightarrow} O_\tau$ between the faces $\delta$ of
$\sigma$, the dual faces\footnote{By abuse of
notation, here $\delta^\bot\cap \sigma^\vee$ is denoted simply by
$\delta^\bot$.} $\tau=\delta^\bot$ of $\sigma^\vee$, and
the $\TT$-orbits $O_\tau$ on~$X$ such that $\dim O_\tau=\dim
\tau=\dim\sigma-\dim\delta$. In particular, the unique $\TT$-fixed point
on $X$ corresponds to the vertex of the cone $\sigma^\vee$, and
the open $\TT$-orbit to the cone $\sigma^\vee$ itself. These
correspondences respect the inclusions: the $\TT$-orbit $O_\mu$
is contained in the orbit closure $\overline{O_\tau}$ if and only if
$\mu\subseteq\tau$, or, equivalently, if and only if $\mu^\bot\supseteq\tau^\bot$
(cf.~\cite[\S 3.2]{CLS}).

For every face $\tau\subseteq\sigma^\vee$, there is a direct sum
decomposition
$$
\kk[X] \ = \kk[\overline{O_\tau}]\oplus I(\overline{O_\tau})\,,
$$
where $
\ \kk[\overline{O_\tau}]=\bigoplus_{m \in
\tau\cap M} \kk\chi^m$ and
\begin{equation}\label{ideal}
I(\overline{O_\tau})=\bigoplus_{m \in (\sigma^\vee\setminus \tau)
\cap M} \kk\chi^m
\end{equation}
is the graded ideal of the subvariety $\overline{O_\tau}$ in $\kk[X]$.

A stabilizer $\TT_{p}=\Stab_\TT(p)$ of any
point $p\in X$ is connected, hence $\TT_{p}\subseteq \TT$ is a
subtorus. Furthermore, $\TT_p\subseteq \TT_q$ if and only if
$\TT.q\subseteq\overline{\TT.p}$, and $\overline{\TT.p}=X^{\TT_p}$;
here $X^{G}$ stands, as usual, for the set
of fixed points of the group~$G$ acting on~$X$.

\medskip

3. {\sl Roots and associate one-parameter unipotent subgroups.}
\begin{definition}[{\normalfont (see~\cite{Dem})}]
A {\it root} of the cone $\sigma$ is a vector $e\in M$ such that
for some index $i$ with $1\le i\le r$, where $r=\card \Xi$, we
have
\begin{equation}\label{old1}
\langle\rho_i, e\rangle \, = \, -1 \quad \text{and} \quad
\langle\rho_j, e\rangle\ge 0 \quad \text{for every} \quad j\ne
i\,.
\end{equation}
\end{definition}

\noindent Let us denote by ${\mathscr R}(\sigma)$ the set of all
roots of the cone~$\sigma$. There is a one-to-one correspondence
$e\stackrel{1\div 1}{\longleftrightarrow} H_e$ between the roots
of~$\sigma$ and the one-parameter unipotent subgroups of $\Aut(X)$
normalized by the torus (see~\cite{Dem} or~\cite{Li}).
Let $\rho_e := \rho_i$. The root $e \in {\mathscr R}(\sigma)$ defines
an LND $\partial_e$ of the $M$-graded algebra $\kk[X]$ given by
\begin{equation}\label{new1}
\partial_e(\chi^m) \ = \ \langle \rho_e,m \rangle \chi^{m+e}\,.
\end{equation}
Its kernel is a (finitely generated) graded subalgebra of $\kk[X]$ (see~\cite{Li}):
\begin{equation}\label{new2}
\ker\partial_e=\bigoplus_{m\in\rho_e^\bot\cap M}\kk\chi^m
\,,
\end{equation}
where $\rho_e^\bot = \{m\in\sigma^\vee
\cap M, \, \langle \rho_e,m\rangle=0\}$ is the
facet (i.e. its codimension one face) of $\sigma^\vee$ orthogonal to $\rho_e$.

\begin{definition}[{(\normalfont see \cite{Fre},~\cite{Li})}]
Two roots~$e$ and~$e'$ with $\rho_e=\rho_{e'}$ are called {\em
equivalent}; we write $e\sim e'$. Two roots~$e$ and~$e'$ are
equivalent if and only if $\ker \partial_e=\ker \partial_{e'}$.
\end{definition}

\begin{remark}\label{eqcl}
Enumerating the ray generators $\Xi=\{\rho_1,\ldots,\rho_r\}$
yields a disjoint partition
$$
{\mathscr R}(\sigma)=\bigcup_{i=1}^r {\mathscr R}_i,
$$
where ${\mathscr R}_i =\{e\in {\mathscr
R}(\sigma)\;\vert\;\rho_e=\rho_i\}$
are nonempty. Indeed,
consider the facet~$\tau_i$ of the cone $\sigma^\vee$, orthogonal to the ray
generator~$\rho_i$. For every~$v$ in the relative interior
$\Int_{\rm rel}(\tau_i)$, the inequalities $\langle\rho_j, v\rangle > 0$ hold for all $j\neq i$. Let $e_0\in M$ be such that
$\langle\rho_i,e_0\rangle =-1$, and let $v_0\in \Int_{\rm
rel}(\tau_i)\cap M$. Putting $e=e_0+kv_0$ for some $k\gg 1$, we obtain
$\langle\rho_j,e\rangle >0$ for all $j\neq i$ and
$\langle\rho_i,e\rangle =-1$. Hence $e\in{\mathscr R}_i$.

\medskip

For example, let $X=\AA^2$ be the affine plane with the standard
$2$-torus action. In this case the cones $\sigma$ and $\sigma^\vee$ coincide with
the first quadrants. The set ${\mathscr R}(\sigma)$ consists of two
equivalence classes
$$
{\mathscr R}_1=\{(x,-1)\,\vert\, x\in\ZZ_{\ge
0}\},\qquad {\mathscr R}_2=\{(-1,y)\,\vert\,
y\in\ZZ_{\ge 0}\}\,.
$$
\end{remark}

4. {\sl One-parameter groups of automorphisms.}
The derivation $\partial_e$ generates a one-parameter
unipotent subgroup $H_e=\lambda_e({\mathbb G}_{\rm a}(\kk)) \subseteq\Aut(X)$,
where $\lambda_e:t\longmapsto \exp(t\partial_e)$. The algebra of
invariants $\kk[X]^{H_e}$ coincides with $\ker\partial_e$.  The
inclusion $\kk[X]^{H_e} \subseteq\kk[X]$ induces a morphism $\pi
\colon X \to Z= \Spec \kk[X]^{H_e}$ whose general fibers are
one-dimensional $H_e$-orbits isomorphic to $\AA^1$ (cf.~\cite[Theorems~2.3 and 3.3]{PV2}). The torus $\TT$ normalizes the
subgroup $H_e$. In particular, $\TT$
stabilizes the fixed point set $X^{H_e}$.

Let $R_e=R_{\rho_e}\subseteq \TT$ be the one-parameter subtorus
corresponding to the vector $\rho_e\in N$. The action  of
$R_e$ on the graded algebra $\kk[X]$ can be given, under a
suitable parametrization $\rho_e:{\mathbb G}_{\rm m}(\kk)\ni t\longmapsto
\rho_e(t)\in R_e$, by
\begin{equation}\label{new3}
t.\chi^m \ = \ t^{\langle \rho_e,m \rangle} \chi^{m},\qquad t\in
{\mathbb G}_{\rm m}(\kk)\,.
\end{equation}
In particular
$\kk[X]^{R_e}=\kk[X]^{H_e}$. Hence the morphism $\pi:X\to Z$
coincides  with the quotient map $X\to X/\!/R_e$. So the
general $H_e$-orbits are the closures of the general $R_e$-orbits. By
Proposition~\ref{pr1} (see below) the latter actually holds for every
one-dimensional $H_e$-orbit.\footnote{If $X$ is a surface then all
fibers of $\pi$ are $R_e$-orbit closures isomorphic to $\AA^1$,
see the parabolic case in~\cite{FZ2}. However, for a toric affine
3-fold~$X$, some degenerate fibers of $\pi$ can be
two-dimensional.}

There is a direct sum decomposition
\begin{equation}\label{new4}
\kk[X]=\kk[X]^{R_e} \ \oplus \
\bigoplus_{m \in \sigma^\vee\cap M \setminus \rho_e^\bot}
\kk \chi^m =\kk[X]^{R_e}\oplus I(D_e)\,,
\end{equation}
where
$D_e:=X^{R_e}\cong Z$.  The divisor $D_e$ coincides with the
attractive set of the $R_e$-action on $X$. So every
one-dimensional $R_e$-orbit has a limit point on $D_e$. The
following simple lemma completes the picture.

\begin{lemma} \label{newl0}
Let $\tau$ be a face of the cone
$\sigma^\vee$, $\Of_\tau$ the corresponding orbit, $\TT_\tau$ the
stabilizer of a point in $\Of_\tau$, and $\Xi_\tau$ the set of ray generators of the
dual face $\tau^\bot\subseteq\sigma$. Then the following holds.
\begin{enumerate}
\item[{\rm a)}] The orbit closure $\overline{\Of_\tau}$ is
stable under $H_e$ if and only if
\begin{equation}\label{eq40}
m+e\in \sigma^\vee\setminus \tau \qquad \forall m\in (\sigma^\vee\setminus \tau)\cap
M : \;\langle
\rho_e,m\rangle>0\,.
\end{equation}
\item[{\rm b)}]  The closure
$\overline{\Of_\tau}$ is $H_{e'}$-stable for any root $e'\sim e$
of the cone $\sigma$ if one of the following equivalent conditions is
fulfilled:
\item[{\rm (i)}] $\rho_e\not\in\Xi_\tau$,
\item[{\rm (ii)}] $\overline{\Of_\tau}\not\subseteq D_e$,
\item[{\rm (iii)}] $R_e\not\subseteq \TT_\tau$.
\end{enumerate}
\end{lemma}

\begin{proof} a) By
virtue of (\ref{new1})
the ideal $I(\overline{\Of_\tau})$ is $\partial_e$-invariant if and only
if
\begin{equation}\label{eq20}
\chi^{m+e} \in I(\overline{\Of_\tau})\qquad \forall \chi^{m} \in
I(\overline{\Of_\tau}) : \;\langle
\rho_e,m\rangle>0\,,
\end{equation}
 which  is equivalent to~(\ref{eq40}) (see (\ref{ideal})). This proves a).

b) For $m\in M$,
$$
m\in \sigma^\vee\setminus\tau\;\iff\; \langle\rho,m\rangle\ge 0
\quad\forall \rho\in\Xi\quad\mbox{and}\quad \exists
\rho\in\Xi_\tau:\;\langle\rho,m\rangle> 0\,.
$$
For any $\rho\neq \rho_e$ we have
$\langle\rho,m+e\rangle\ge\langle\rho,m\rangle$. Hence
(i) implies~(\ref{eq40}).

We have $\overline{\Of_\tau}=X^{\TT_\tau}$ and
$D_e=X^{R_e}=\overline{\Of_{\rho_e^\bot}}$, where
$\rho_e^\bot=(\RR_+\rho_e)^\bot$. So the equivalence
(i)$\iff$(ii)$\iff$(iii) is due to the
Orbit-Cone correspondence.
\end{proof}

\begin{remark}\label{remk1}
Consequently, $\rho_e\in \Xi_\tau$ if
$\overline{\Of_\tau}$ is not $H_e$-stable. In general, the converse is not
true. For instance, let $X=\AA^2$ be the plane  with
the standard 2-dimensional torus action. Here $\Xi=\{(1,0),(0,1)\}$. Let
$\tau=\{ 0\}$ and $e=(0,-1)$, $e'=(a,-1)\sim e$, where $a>0$
and $\rho_e=(0,1)$. Then (\ref{eq40}) holds for $H_{e'}$ but not
for $H_e$. Hence $\rho_e\in \Xi_\tau$, and the $\TT$-fixed point
$\overline{\Of_\tau}=\{(0,0)\}$ is $H_{e'}$-stable but not
$H_{e}$-stable.
One can also construct an example with $\dim X=4$ such that the
closure $\overline{\Of_\tau}$ is $H_{e'}$-stable for every root
$e'\sim e$, whereas the equivalent conditions (i)--(iii) are not
fulfilled.
\end{remark}

For the proof of infinite transitivity we need somewhat more
precise information concerning the actions of one-parameter groups of
automorphisms on toric varieties; see Proposition~\ref{pr1} and Lemmas~\ref{l1},~\ref{l2} below.

\begin{proposition}\label{pr1}
Given a root $e\in {\mathscr R}(\sigma)$.
Let, as before, $H_e\subseteq\SAut (X)$ be the associate
one-parameter unipotent subgroup. Then the following holds.
\begin{enumerate}
\item[{\rm a)}]
For every point $x\in X \setminus X^{H_e}$, the orbit $H_e. x$
meets exactly two $\TT$-orbits $\Of_1$ and $\Of_2$ on $X$, where
$\dim \Of_1 = 1+\dim \Of_2$.
\item[{\rm b)}]
The intersection $\Of_2\cap H_e. x$ consists of a single point, while
$$
\Of_1\cap H_e. x \ = \ R_e. y\qquad \forall y \in \Of_1\cap H_e.
x\,.
$$
\end{enumerate}
\end{proposition}

\begin{proof} a)
The number of $\TT$-orbits in $X$ being finite, there exists a
$\TT$-orbit $\Of_1$ such that the intersection $\Of_1 \cap H_e . x$ is a
Zariski open subset of $H_e . x$. So $H_e . x\subseteq
\overline{\Of_1}$. There is also another $\TT$-orbit $\Of_2$ that
meets $H_e . x$. Indeed, otherwise $H_e . x\cong \AA^1$ would be
contained in a single $\TT$-orbit $\Of_1$. However, this is
impossible because the regular invertible functions separate points on $\Of_1$, and, consequently, either they do on $H_e . x$. Since $\Of_2$ intersects
$\overline{\Of_1}$, we have $\Of_2\subseteq \overline{\Of_1}$ and
so $\dim \Of_2<\dim \Of_1$.

We know that the torus $\TT$ normalizes the unipotent subgroup $H_e$. Hence the
elements of $\TT$ send the $H_e$-orbits into $H_e$-orbits. In
particular, for every point $q\in H_e . x$ the stabilizer $\TT_{q}$
preserves the orbit $H_e . x$. For all $q\in\Of_1 \cap H_e . x$
the stabilizer is the same. Since $H_e . x\subseteq
\overline{\Of_1}=X^{\TT_q}$, this stabilizer acts trivially on $H_e . x$.
Thus $\TT_r\supseteq \TT_q$ for any point $r\in H_e . x$, and
$\TT_r=\TT_q$ if and only if $r\in \Of_1 \cap H_e . x$.

Fix a point $p\in\Of_2 \cap H_e . x$. If $\TT_p\subseteq \TT_q$, then
$\TT_p=\TT_q$, and so $\dim \Of_2=\dim \Of_1$, a contradiction.
Consequently, the stabilizer $\TT_p$ acts on $H_e . x$ with two
orbits, i.e., $H_e. x = \TT_{p}\, . q \cup \{p\}$, where $q\in H_e .
x\setminus\{p\}$. From the exact sequence
$$
1\to \TT_q\to \TT_p\to {\mathbb G}_{\rm m}(\kk)\to 1
$$
we deduce that $\dim \TT_p=1+\dim \TT_q$. Finally
$H_e. x \subseteq \Of_1 \cup \Of_2$, and $\dim \Of_1 = 1+ \dim
\Of_2$, as stated in a).

b) We may assume that $\Of_1=\TT . x$. Since $H_e . x
\subseteq\overline{\Of_1}$ and the torus $\TT$ normalizes the
subgroup $H_e$, we have $H_e(\Of_1)\subseteq\overline{\Of_1}$. Thus
$\overline{\Of_1}$ is $H_e$-stable. On the other hand, since
$H_e.p=H_e.x\not\subseteq \overline{\Of_2}$, the closure
$\overline{\Of_2}$ is not $H_e$-stable. In particular, by Lemma~\ref{newl0}
$\rho_e\in\Xi_{\tau_2}$, where $\tau_i$ is the face of
the cone $\sigma^\vee$ which corresponds to $\Of_i=\Of_{\tau_i}$, $i=1,2$.
Hence by the same lemma $R_e\subseteq \TT_{\tau_2}=\TT_p$. Let us show
that $R_e\not\subseteq \TT_{\tau_1}=\TT_q$, where $q\in H_e.x\setminus\{p\}$.
Applying again Lemma~\ref{newl0}, we obtain that (\ref{eq40}) holds
for $\tau=\tau_1$ but not for $\tau=\tau_2$. Since
$\tau_2\subseteq\tau_1$, this implies that (\ref{eq40}) does not
hold for some $m\in\tau_1\setminus\tau_2$. The latter is
possible only if $\rho_e\in \tau_2^\bot\setminus\tau_1^\bot$. Thus by
Lemma~\ref{newl0} $R_e\not\subseteq \TT_{\tau_1}=\TT_q$.
Finally, the one-dimensional orbit $\TT_{p}. q$ coincides with
$R_e.q$.
\end{proof}

\begin{definition}\label{con}
We say that a pair of $\TT$-orbits $(\Of_1, \Of_2)$ in~$X$ is {\it
$H_e$-connected} if $H_e . x\subseteq \Of_1\cup \Of_2$ for some
$x\in X\setminus X^{H_e}$. By Proposition~\ref{pr1}
$\Of_2\subseteq\overline{\Of_1}$ for such a
pair (up to a permutation), and
$\dim\Of_1=1+\dim\Of_2$. Clearly, we can choose a point~$x$ on the
orbit $\Of_2$, as above. Since the torus normalizes the
subgroup $H_e$, any point of $\Of_2$ can actually serve as such a
point~$x$.
\end{definition}

\begin{example}\label{re-new}
Given a root $e\in {\mathscr R}(\sigma)$,
the derivation $\partial_e$ as in (\ref{new1}) extends to an LND
of a bigger graded algebra
$$
A(\rho_e)=\bigoplus_{m\in M, \langle\rho_e,m\rangle\ge 0}\kk\chi^m\,.
$$
Indeed, letting
$k=\langle\rho_e,m\rangle\ge 0$, we yield
$\langle\rho_e,m+ke\rangle=0$, and so $\partial_e^k
(\chi^m)\in\ker\partial_e$. This provides a $\TT$- and
$H_e$-stable open subset
$$
U=\Spec A(\rho_e)\cong (\AA^1_*)^{n-1}\times\AA^1\subseteq X\,,
$$
where $n=\dim X$, $\AA^1_*=\Spec \kk[t,t^{-1}]$, $\AA^1=\Spec \kk[u]$, here
$u=\chi^{-e}$, and $H_e$ acts along the second factor by the shifts.
The only $\TT$-orbits in $U$ are the open orbit
$\Of_1=\{u\neq 0\}$, which corresponds to the vertex of $\sigma$,
and the codimension one orbit $\Of_2=\{u=0\}$, which corresponds
to the ray $\QQ\rho_e$. It is easy to see that the
pair $(\Of_1,\Of_2)$ is $H_e$-connected.
\end{example}

From Proposition~\ref{pr1} and its proof we deduce the following
criterion of $H_e$-connectedness.

\begin{lemma} \label{l1}
Let $(\Of_1,\Of_2)$ be a pair of $\TT$-orbits on~$X$ with
$\Of_2\subseteq\overline{\Of_1}$, where
$\Of_i=\Of_{\sigma_i^\bot}$ for a face $\sigma_i$ of the cone $\sigma$,
$i=1,2$. Given a root $e\in {\mathscr R}(\sigma)$, the pair
$(\Of_1,\Of_2)$ is $H_e$-connected if and only if
$e\vert_{\sigma_2}\le 0$ and $\sigma_1$ is a facet of $\sigma_2$
given by the equation $\langle v,e\rangle=0$.
\end{lemma}

\begin{proof}
In course of the proof of Proposition~\ref{pr1},\,b) we established
that the pair $(\Of_1,\Of_2)$ is $H_e$-connected if and only if the closure
$\overline{\Of_1}$ is $H_e$-invariant, $\overline{\Of_2}$ is not,
and $\dim \Of_1=1+\dim \Of_2$. Moreover, if the pair $(\Of_1,\Of_2)$ is
$H_e$-connected then $\sigma_2^\bot$ is a facet of the cone $\sigma_1^\bot$
(and so $\sigma_1$ is a facet of $\sigma_2$), and there exists
$m_0\in \sigma_1^\bot\setminus \sigma_2^\bot$ such that
$$
\langle \rho_e,m_0\rangle>0, \qquad m_0+e\in\sigma_2^\bot\,.
$$
Since $\langle \rho_i,e\rangle\ge 0$ for all $\rho_i\neq\rho_e$, we obtain that
$\sigma_2=\cone (\sigma_1,\rho_e)$. We have also
$e\vert_{\sigma_1}=0$ because $e=m_0+e-m_0\in\span \sigma_1^\bot$.
Thus $e\vert_{\sigma_2}\le 0$ and $\sigma_1$ is given in
$\sigma_2$ by the equation $\langle v,e\rangle=0$.

Conversely, assume that $e\vert_{\sigma_2}\le 0$ and $\sigma_1$ is
given in $\sigma_2$ by the equation $\langle v,e\rangle=0$. Then for
any $m\in\sigma^\vee\setminus \sigma_1^\bot$ such that $\langle
\rho_e,m\rangle>0$ we have $m+e\not\in \sigma_1^\bot$. Indeed,
$e\vert_{\sigma_1}=0$ and so $e\in\sigma_1^\bot$. Thus the condition
(\ref{eq40}) holds for $\sigma_1^\bot$. Furthermore, $\langle
\rho_e,m'\rangle>0$ for any $m'\in  \sigma_1^\bot\setminus
\sigma_2^\bot$. It follows that
$$
m_0:=m'+(\langle \rho_e,m'\rangle-1)\cdot e\in
\sigma_1^\bot\setminus \sigma_2^\bot, \qquad
m_0+e\in\sigma_2^\bot\,.
$$
Indeed, $\langle \rho_e,m_0\rangle=1$
and $\langle \rho_e,m_0+e\rangle =0$, while $\langle \rho_i,m_0\rangle\ge 0$
and $\langle \rho_i,m_0+e\rangle\ge 0$ for
every $\rho_i\neq\rho_e$. Therefore (\ref{eq40}) is fulfilled for
$\sigma_1^\bot$ but not for $\sigma_2^\bot$. Consequently, the
pair $(\Of_1,\Of_2)$ is $H_e$-connected.
\end{proof}

\begin{remark}\label{orclo}
Given a one-parameter subgroup $R\subseteq\TT$ and a point $x\in
X\setminus X^R$, the orbit closure $\overline{R.x}$ coincides with
an $H_e$-orbit if and only if $\overline{R.x}$ is covered by a
pair of $H_e$-connected $\TT$-orbits. For instance, for $X=\AA^2$
with the standard 2-dimensional torus action and $R\subseteq\TT$ being the
subgroup of scalar matrices, the latter condition holds only for
the points $x\ne 0$ on one of the coordinate axes.
\end{remark}

\begin{lemma} \label{l2}
For any point\,\footnote{Recall that $\Of_{\sigma^\vee}$ is the open $\TT$-orbit in
$X$.} $x\in X_{\reg} \setminus \Of_{\sigma^\vee}$
there is
a root $e\in {\mathscr R}(\sigma)$ such that
$$
\dim \TT.y \ > \ \dim \TT. x
$$
for a general point $y\in H_e.x$. In particular, the pair
$(\TT.y,\TT. x)$ is $H_e$-connected.
\end{lemma}

\begin{proof}
Since $x\not\in \Of_{\sigma^\vee}$, by the Orbit-Cone
Correspondence there exists a proper face, say,
$\sigma_2\subseteq\sigma$ such that $\TT.x=\Of_{\sigma_2^\bot}$.
The point $x\in X$ being regular, the ray generators, say,
$\rho_1,\ldots,\rho_s$ of the cone $\sigma_2$ form a base of a primitive
sublattice $N'\subseteq N$ (see~\cite[\S 2.1]{Fu}). Let
$\sigma_1$ be the facet of $\sigma_2$ spanned by
$\rho_2,\ldots,\rho_s$. Again by the Orbit-Cone Correspondence,
$\Of_{\sigma_2^\bot}\subseteq\overline{\Of_{\sigma_1^\bot}}$ and
$\dim \Of_{\sigma_1^\bot} = 1+\dim \Of_{\sigma_2^\bot}$. Let us
show that the pair $(\Of_{\sigma_1^\bot}, \Of_{\sigma_2^\bot})$ is
$H_e$-connected for some root $e\in {\mathscr R}(\sigma)$
satisfying the assumptions of Lemma~\ref{l1}.

Choosing a $\sigma$-supporting hyperplane $L \subseteq N_{\QQ}$ such
that
$\sigma_2 = \sigma \cap L$, we obtain a splitting $N=N'\oplus N'' \oplus N'''$,
where $N \cap L=N'\oplus N''$ and $N'''\cong \ZZ$. Consider
a linear form $e_1$ on $N'$ defined by
$$
\langle \rho_1,e_1\rangle\, =\, -1 , \qquad \langle
\rho_2,e_1 \rangle\ =\ \ldots\ = \ \langle \rho_s,e_1\rangle\ =\
0\,.
$$
Let $e_2$ be a non-zero linear form on $N'''$. Extending~$e_1$ and~$e_2$
to the whole lattice~$N$ by zero on the complementary
sublattices, we obtain a linear form $e = e_1+e_2$ on~$N$.
Multiplying $e_2$ by a suitable integer, we can achieve that
$\langle \rho_j,e\rangle>0$ for every $\rho_j\notin\sigma_2$. Then~$e$
is a root of the cone $\sigma$ such that $\rho_e=\rho_1$ and
the condition of Lemma~\ref{l1} holds for~$e$. By Lemma~\ref{l1}
the pair $(\Of_{\sigma_1^\bot}, \Of_{\sigma_2^\bot})$
is $H_e$-connected, as claimed. Since $\TT.x=\Of_{\sigma_2^\bot}$
and the torus $\TT$ normalizes the group $H_e$,
the desired conclusion follows from Proposition~\ref{pr1} and
the observation in Definition~\ref{con}.
\end{proof}

Now we possess all the necessary notions to prove infinite transitivity in Theorem~\ref{T2}. The proof consists of several steps; see
Lemmas~\ref{l3}-\ref{flex} below.

\begin{lemma} \label{l3}
For any collection of $m$ distinct points $Q_1,\ldots,Q_m\in X_{\reg}$ there exists an
automorphism $\phi \in \SAut(X)$ such that the images
$\phi(Q_1),\ldots,\phi(Q_m)$ are contained in the open $\TT$-orbit.
\end{lemma}

\begin{proof}
Let
$$
d(Q_1,\ldots,Q_m) = \dim \TT. Q_1 + \ldots + \dim \TT. Q_m
$$
and assume that $\dim \TT. Q_i < \dim X$ for some~$i$. By
Lemma~\ref{l2} there exists a root $e\in {\mathscr R}(\sigma)$ such that
$\dim \TT.P_i > \dim \TT. Q_i$ for a general point $P_i \in H_e.Q_i$. Fix an isomorphism
$\lambda_e:{\mathbb G}_{\rm a}(\kk)\stackrel{\cong}{\longrightarrow} H_e$. There is
a finite set of values $t\in {\mathbb G}_{\rm a}(\kk)$ such that $\dim \TT.(\lambda_e(t).Q_j) \
< \ \dim \TT. Q_j$ for some $j\ne i$. Thus for a general $t\in
{\mathbb G}_{\rm a}(\kk)$ we have
$$
d(\lambda_e(t).Q_1,\ldots,\lambda_e(t).Q_m)> d(Q_1,\ldots,Q_m)\,.
$$
Applying recursion, we get the result.
\end{proof}

From now on we assume that the points $Q_1,\ldots,Q_m$ are contained in the
open $\TT$-orbit $\TT. x_0$. We fix a maximal subset of
linearly independent ray generators
$\{\rho_1,\ldots,\rho_n\}=:\Xi^{(0)} \subseteq \Xi$, where $n=\dim X$.
For every $i=1,\ldots,n$  we choose
an isomorphism $\rho_i:{\mathbb G}_{\rm m}(\kk)\stackrel{\cong}{\longrightarrow}
R_{\rho_i}$, denoted by the same letter as the ray generator.
Recall that for a root $e\in {\mathscr R}(\sigma)$ the inclusion
$\kk[X]^{H_e} \subseteq\kk[X]$ induces a morphism $\tau \colon X
\to Z$, where $Z = \Spec \kk[X]^{H_e}$.

\begin{lemma} \label{l4}
Given $\rho_i\in\Xi^{(0)}$ and a root $e\in {\mathscr R}(\sigma)$ such that $\rho_e=\rho_i$. Then for every finite set
$\Lf_0,\ldots,\Lf_k$ of pairwise distinct $R_e$-orbits in
$\TT.x_0$ there exists a regular invariant $q\in\kk[X]^{H_e}$
which is identically equal to $1$ on $\Lf_0$ and vanishes on
$\Lf_1,\ldots,\Lf_k$.
\end{lemma}

\begin{proof}
The quotient morphism $\tau \colon X \to Z$ separates typical
$H_e$-orbits (see \cite[Theorems~2.3 and 3.3]{PV2}). Since the
torus $\TT$ normalizes the subgroup $H_e$, there exists a $\TT$-action
on~$Z$ such that the morphism $\tau$ is $\TT$-equivariant. In
particular, for every $x\in \TT. x_0$ the fiber of $\tau$ passing
through~$x$ is an $H_e$-orbit. According to Proposition~\ref{pr1}, the
$R_e$-orbits $\Lf_0,\dots,\Lf_k$ are intersections of the
corresponding $H_e$-orbits with the open orbit $\TT. x_0$. Thus
for every $j=1,\ldots,k$ there exists an invariant $q_j\in
\kk[X]^{H_e}$ which vanishes on $\Lf_j$ and restricts to the constant function~$1$
on $\Lf_0$. It is easy to see that the product
$q=q_1\cdot\ldots\cdot q_k\in \kk[X]^{H_e}$ has the desired
properties.
\end{proof}

In the notation of Lemma~\ref{l4}, let
$\Stab_{\Lf_1,\ldots,\Lf_k}(\Lf_0)\subseteq \SAut(X)$ denote the
subgroup of all transformations that fix the orbits
$\Lf_1,\ldots,\Lf_k$ pointwise and stabilize the closure $\overline{\Lf_0}$
in~$X$.

\begin{lemma} \label{l5}
There exists a one-parameter unipotent subgroup
$H\subseteq\Stab_{\Lf_1,\ldots,\Lf_k}(\Lf_0)$ which acts
transitively on $\overline{\Lf_0}$.
\end{lemma}

\begin{proof}
As before, let $e\in {\mathscr R}(\sigma)$ be a root  with
$\rho_e=\rho_i$, and~$q$ be a regular $H_e$-invariant as in
Lemma~\ref{l4}. The LND $q\partial_e\in\Der\kk[X]$ defines a
one-parameter unipotent subgroup
$H\subseteq\Stab_{\Lf_1,\ldots,\Lf_k}(\Lf_0)$. Clearly, the
restriction $H\vert_{\Lf_0}=H_e\vert_{\Lf_0}$ acts by shifts on
$\overline{\Lf_0}\cong\AA^1$, and hence this action is transitive.
\end{proof}

In the remaining part of the proof of Theorem~\ref{T2} we use the
following notation. Let
$\Xi^{(0)}=\{\rho_1,\ldots,\rho_n\}$ be a basis in $N_{\QQ}$ formed by ray
generators. We consider the homomorphism $\theta \colon
{\mathbb G}_{\rm m}(\kk)^n \, \to \, \TT$ of the standard $n$-torus to $\TT$
given by
\begin{equation}\label{one}
\theta: (t_1,\ldots,t_n) \, \longmapsto \, (\rho_1(t_1)
\cdot\ldots\cdot \rho_n(t_n))\,.
\end{equation}
It is easily seen that $\theta$ is surjective and its kernel
$\Theta=\ker (\theta)$ is a finite subgroup in ${\mathbb G}_{\rm m}(\kk)^n$.
We consider as well the induced surjective morphism from
${\mathbb G}_{\rm m}(\kk)^n$ to the open orbit $\TT. x_0$. In particular,
for~$m$ distinct points $Q_1,\ldots,Q_m\in \TT.x_0$ we can write
\begin{equation}\label{mark}
Q_j \ = \ \theta(t_{1,j},\ldots,t_{n,j}). x_0, \qquad
j=1,\ldots,m\,,
\end{equation}
where the point $(t_{1,j},\ldots,t_{n,j})\in {\mathbb G}_{\rm m}(\kk)^n$ is
determined by~$Q_j$ up to the diagonal action of $\Theta$ on
${\mathbb G}_{\rm m}(\kk)^n$:
\begin{equation}\label{two}
\vartheta.(t_1,\ldots,t_n)
=(\vartheta_1t_1,\ldots,\vartheta_nt_n), \qquad\mbox{where}\quad
\vartheta=(\vartheta_1,\ldots,\vartheta_n)\in\Theta\,.
\end{equation}
Letting $\kappa=\ord \Theta$, by the Lagrange Theorem we have
$\vartheta_i^\kappa=1$ for all $i=1,\ldots,n$.

We fix a standard collection of~$m$ points in $\TT.x_0$:
\begin{equation}\label{mark1}
Q_{j0} \, = \, \theta(j,\ldots,j).x_0,\qquad
j=1,\ldots,m\,.
\end{equation}
By our assumptions $\Char (\kk)=0$
and $\overline{\kk} =\kk$, hence these points are distinct. It
remains to find an automorphism $\varphi \in \SAut(X)$ such
that $\varphi(Q_j) \, = \, Q_{j0}$ for every $j=1,\ldots,m$. To
this end we use Lemmas~\ref{l6} and~\ref{l7} (see below).

We say that $t,t'\in {\mathbb G}_{\rm m}(\kk)$ are {\em $\kappa$-equivalent}
if $t'=\varepsilon t$ for some $\kappa$th
root of unity $\varepsilon\in {\mathbb G}_{\rm m}(\kk)$.

\begin{lemma}\label{l6}
\begin{enumerate}
\item[{\rm a)}] For any pairwise distinct
elements $t_1,\ldots,t_n\in {\mathbb G}_{\rm m}(\kk)$, the set of values $a\in\kk$
such that $t_i+a$ and $t_j+a$ are $\kappa$-equivalent for some
$i\neq j$ is finite.
\item[{\rm b)}] Fix $s\in\{1,\ldots,n\}$. If the
points $Q_i$ and $Q_j$ belong to the same $R_s$-orbit, then their
$r$th components $t_{i,r}$ and $t_{j,r}$ are $\kappa$-equivalent
for every $r\neq s$.
\item[{\rm c)}] Suppose that the points $Q_{j_1},\ldots, Q_{j_l}$
belong to the same $R_s$-orbit $\Lf^{(s)}$. Then their images under a
general shift on the line $\overline{\Lf^{(s)}}\cong\AA^1$ belong
to distinct $R_r$-orbits for every $r\neq s$.
\end{enumerate}
\end{lemma}

\begin{proof}
a) Given a $\kappa$th root of unity $\varepsilon$,
the linear equation
$$
t_i+a=\varepsilon (t_j+a)
$$
is satisfied for at most one value of~$a$. It implies~a).

The assertion of~b) holds since the $R_s$-action on~$X$, lifted
via~(\ref{one}), affects only the component $t_{i,s}$ of the point~$Q_i$ in
(\ref{mark}), while the $\Theta$-action on ${\mathbb G}_{\rm m}(\kk)^n$
replaces the component $t_{i,r}$ ($r\neq s$) by a
$\kappa$-equivalent one.

Now c) follows immediately from~a) and~b). Indeed, for $i\neq j$ the
intersection of any $R_i$- and $R_j$-orbits is at most finite.
\end{proof}

\begin{lemma}\label{l7}
In the notation as above there exists
$\psi\in\SAut (X)$ such that the points
$\psi(Q_1),\ldots,\psi(Q_m)$ belong to different $R_1$-orbits in
$\TT.x_0$.
\end{lemma}

\begin{proof}
By our assumption $n\ge 2$,
so there is an $R_2$-action on~$X$. Let
$\Lf_0^{(2)},\ldots,\Lf_k^{(2)}$ be the distinct $R_2$-orbits
passing through the points $Q_1,\ldots,Q_m$ so that this
collection splits into $k+1$ disjoint pieces. We may assume that
the piece on $\Lf_0^{(2)}$ is $Q_1,\ldots,Q_l$. Applying Lemma~\ref{l5}
to $\rho_e=\rho_2$, we find a one-parameter
unipotent subgroup
$H\subseteq\Stab_{\Lf_1^{(2)},\ldots,\Lf_k^{(2)}}(\Lf_0^{(2)})$
acting by shifts on $\overline{\Lf_0^{(2)}}\cong\AA^1$. By Lemma~\ref{l6}
the images of $Q_1,\ldots,Q_l$ under a general shift
lie in different $R_r$-orbits for every $r\neq 2$, while all the
other points $Q_j$, $j>l$, remain fixed. Applying the same
procedure subsequently to the other pieces, we finally obtain an automorphism
$\psi\in \SAut (X)$ such that the points
$\psi(Q_1),\ldots,\psi(Q_m)$ belong to different $R_r$-orbits for
every $r\neq 2$.
\end{proof}

{\sc Proof of Theorem~\ref{T2}.} Let us start with the infinite transitivity.
By virtue of Lemma~\ref{l7} we may assume  that
the orbits $\Lf_j^{(1)}=R_1.Q_j$, $j=1,\ldots,m$, are all
distinct. By Lemma~\ref{l5} we can change the component $t_{1,j}$
of a point $Q_j$ arbitrarily while fixing the other components and
as well the other points of our collection. Thus we can achieve
that $t_{1,j}=j$ for all $j=1,\ldots,m$. This guarantees that for
any $l\ge 2$ the orbits $R_l.Q_1,\ldots,R_l.Q_m$ are pairwise
distinct. Applying Lemma~\ref{l5} again to every $R_l$-orbit for
$l=2,\ldots,n$, we can reach the standard collection
$Q^{(0)}_1,\ldots,Q_m^{(0)}$ as in (\ref{mark1}) with $t_{l,j}=j$
for all $j=1,\ldots,m,\,\,l=1,\ldots,n$. This proves the infinite
transitivity statement in Theorem~\ref{T2}.

The proof of flexibility is in the next lemma.

\begin{lemma}\label{flex}
Every non-degenerate affine toric variety~$X$ is flexible.
\end{lemma}

\begin{proof}
If $\dim X=1$, then $X\cong \AA^1$, and
the assertion is evidently true. Further we suppose that $\dim X\ge
2$. We already know that the group $\SAut (X)$ acts (infinitely)
transitively on $X_{\reg}$. Hence it is enough to find just one
flexible point on $X_{\reg}$. Let us show that a point $x_0$ in
the open $\TT$-orbit is flexible. Consider the action of the
standard torus ${\mathbb G}_{\rm m}(\kk)^n$ on $X$ induced by the $\TT$-action
on $X$ via~(\ref{one}). The stabilizer $\Stab (x_0)\subseteq
{\mathbb G}_{\rm m}(\kk)^n$ being finite, the tangent map
$T_g {\mathbb G}_{\rm m}(\kk)^n\to T_{x_0}X$ is surjective at each point $g\in \Stab (x_0)$.
Hence the tangent vectors at $x_0$ to the orbits
$R_i.x_0$, $i=1,\ldots,n$, span the tangent space $T_{x_0}X$. By
Remark~\ref{eqcl}, for every $i=1,\ldots,n$ there exists a root
$e_i\in{\mathscr R}(\sigma)$ such that $\rho_i=\rho_{e_i}$. Since~$x_0$
cannot be fixed by the one-parameter unipotent subgroup
$H_{e_i}$, it follows from Proposition~\ref{pr1} that $H_{e_i} .
x_0=\overline{R_i.x_0}$. If one parameterizes properly these two orbits,
their velocity vectors at~$x_0$ coincide. Therefore, $T_{x_0}X$ is
spanned as well by the tangent vectors of the orbits $H_{e_i} . x_0$,
$i=1,\ldots,n$, which means that the point $x_0$ is flexible
on~$X$.
\end{proof}

Now the proof of Theorem~\ref{T2} is completed.

\begin{example} \label{rm1}
Consider a singular affine toric surface
$X_{d,e}= \AA^2/G_d$, where $d$ and $e$ are coprime integers with
$0< e<d$, and  $G_d$ is the cyclic group generated by a primitive
$d$th root of unity $\zeta$, which acts on the plane $\AA^2$ via
$\zeta. (x,y) = (\zeta x, \zeta^e y)$. It is well known (see \cite{FZ2}, \cite{Gi}, and~\cite{Po})
that for $e\ge 2$ the smooth
locus $(X_{d,e})_{\reg}=X_{d,e}\setminus\{0\}$ is not isomorphic
to a homogeneous space of any affine algebraic group. However
$X_{d,e}\setminus\{0\}$ is homogeneous under the action of the
infinite dimensional group $\SAut (X)$.
\end{example}

\section{Affine suspensions}

In this section we prove part~3 of Theorem~\ref{main0}. Let us first recall some necessary notions.

\begin{definition}\label{susp}
Let $X^{(0)}$
be an affine variety. By a {\it cylinder} over $X^{(0)}$ we mean
the product $X^{(0)}\times \AA^1$. Given a nonconstant regular
function $f_1\in\kk[X^{(0)}]$, we define a new affine variety
$$
X^{(1)}=\SP(X^{(0)},f_1):=\{f_1(x)-uv=0\}\subseteq X^{(0)}\times
\AA^2,
$$
called a {\it suspension} over $X^{(0)}$. By recursion, for any
$l\in\NN$ we obtain the iterated suspension
$X^{(l)}=\SP(X^{(l-1)},f_{l})$.
\end{definition}

For instance, starting with $X^{(0)}=\AA^k$, we arrive at the $l$th
suspension $X^{(l)}$ given in the affine space
$$
\AA^{k+2l}=\Spec \kk[x_1,\ldots,x_k,u_1,v_1,\ldots,u_l,v_l]
$$
by the system of equations
\begin{equation}\label{system}
\left\{
\begin{aligned}
u_1v_1&-f_1(x_1,x_2,\dots,x_k)=0\\
u_2v_2&-f_2(x_1,x_2,\dots,x_k, u_1,v_1)=0\\
&\ldots \\
u_lv_l&-f_l(x_1,x_2,\dots,x_k,u_1,v_1,u_2,v_2,
\dots,u_{l-1},v_{l-1})=0,\\
\end{aligned}
\right.
\end{equation}
where  for every $i=1,2,\ldots,l$ the polynomial
$f_i\in\kk[x_1,\ldots,x_k,u_1,v_1,\ldots, u_{i-1},v_{i-1}]$ is
non-constant modulo the ideal
$(u_1v_1-f_1,\ldots,u_{i-1}v_{i-1}-f_{i-1})$.

\smallskip

We prove part~3 of Theorem~\ref{main0} separately for suspensions over a line and over a variety of dimension at least~2. For suspensions over a line, the
assertion remains true over an arbitrary field of characteristic
zero, under an additional restriction on the function $f=f_1$.

\begin{theorem}\label{aa}
Let $\kk$ be a field of characteristic  zero and let $f\in\kk[x]$ be a polynomial such that $f(\kk)=\kk$. Consider a surface $X\subseteq
\AA_\kk^3$ given by the equation $f(x)-uv=0\,.$ Then~$X$ is flexible and
the special automorphism group $\SAut(X)$ acts $m$-transitively on
$X_{\reg}$ for every $m\in \NN$.
\end{theorem}

In higher dimensions, part 3 of Theorem~\ref{main0} can be
restated as follows.

\begin{theorem}\label{main}
Let $\kk$ be an
algebraically closed field of characteristic  zero. Let $X^{(0)}$
be a flexible affine variety. Assume that either
$X^{(0)}\cong\AA^1$, or $\dim X^{(0)}\geqslant 2$ and the special
automorphism group $\SAut(X^{(0)})$ acts $m$-transitively on
$X^{(0)}_{\reg}$ for every $m\in \NN$. Then every iterated
suspension $X^{(l)}$ over $X^{(0)}$, $l\ge 1$, is flexible and the
special automorphism group $\SAut(X^{(l)})$ acts $m$-transitively
on $X^{(l)}_{\reg}$ for every $m\in \NN$.
\end{theorem}

Since the assumptions of the theorem are fulfilled for the affine
space $X^{(0)}=\AA^k$, $k\ge 1$, we can conclude that for every
$k,l\ge 1$ the affine variety $X^{(l)}\subseteq \AA^{k+2l}$
defined by (\ref{system}) is flexible, and the group
$\SAut(X^{(l)})$ acts infinitely transitively on $X^{(l)}_{\reg}$.

The proof of Theorem~\ref{main}, with minor changes, works also
for real algebraic varieties and leads to the following result.

\begin{theorem}\label{main10}
Let $X^{(0)}$ be a flexible
real affine algebraic variety. Suppose that the smooth locus
$X_{\reg}^{(0)}$ is connected and the special automorphism group
$\SAut(X^{(0)})$ acts $m$-transitively on $X^{(0)}_{\reg}$ for
every $m\in \NN$. Consider the iterated suspensions
$X^{(i)}=\SP(X^{(i-1)},f_{i})$, where the functions
$f_i\in\RR[X^{(i-1)}]$ satisfy the conditions
$f_i(X^{(i-1)}_{\reg})=\RR$, $i=1,\ldots,l$. Then for every
$i=1,\ldots,l$ the variety $X^{(i)}$ is flexible and the special
automorphism group $\SAut(X^{(i)})$ acts $m$-transitively on
$X^{(i)}_{\reg}$ for every $m\in \NN$.
\end{theorem}

The infinite transitivity in Theorems~\ref{aa}--\ref{main10} is proved in Subsections 3.1--3.3,
respectively. The flexibility in all the three cases is established in
Subsection~3.4.


\subsection{Suspensions over a line}
\begin{proof}[Proof of infinite transitivity in Theorem~\ref{aa}] (The
proof is elementary and is based on some explicit formulae from~\cite{ML1}.)
We may assume that $d=\deg f\ge 2$.
According to\footnote{Cf.\ also~\cite{Da}, \cite{FKZ3}, and~\cite{ML2}.}~\cite{ML1}, in our case the special automorphism group
$\SAut(X)$ contains the abelian subgroups $G_u$ and $G_v$
generated, respectively, by the one-parameter unipotent subgroups
\begin{equation} \label{type-v}
H_u(q): (x, u, v) \mapsto \left(x+tq(u), u,
v+\frac{f(x+tq(u))-f(x)}{u}\right),
\end{equation}
\begin{equation} \label{type-u}
H_v(q): (x, u, v) \mapsto  \left(x+tq(v),
u+\frac{f(x+tq(v))-f(x)}{v}, v\right),
\end{equation}
where $q(z) \in \kk[z]$, $q(0)=0$, and $t\in\kk$. Clearly,
$u\in\kk[X]^{G_u}$ and $v\in\kk[X]^{G_v}$. We claim that the
subgroup $G=\langle G_u,G_v\rangle\subseteq\SAut(X)$ acts
$m$-transitively on $X_{\reg}$ for every $m\in \NN$. So given an
$m$-tuple of pairwise distinct points of~$X_{\reg}$
$$
Q_1=(x_1,u_1,v_1),\quad \ldots,\quad Q_m=(x_m,u_m,v_m),
$$
our aim is to find an automorphism $\phi\in G$ which sends this $m$-tuple
to a standard $m$-tuple
$$
Q_i^{(0)}=(x^{(0)}_i,u^{(0)}_i,v^{(0)}_i),\qquad i=1,\ldots, m\,,
$$
chosen in such a way that all $v_i^{(0)}$ are nonzero and
distinct.

{\sl Step} 1. Acting with the subgroup $G_u$, we can replace the original
$m$-tuple by another one such that $v_i\neq 0$ for all
$i=1,\ldots,m$. Indeed, since $f^{(d)}(x)$ is a nonzero constant, and since
either $f'(x)\neq 0$ or $u\ne 0$, the polynomial
$$
\frac{f(x+tu)-f(x)}{u} = \frac{f'(x)}{1!}t+\ldots+
\frac{f^{(d)}(x)}{d!}u^{d-1}t^d\in \kk[x,u][t]\,,
$$
is non-constant. Since the point $Q_s\in X$ is
smooth, the equalities $u_s=0,\,v_s=0,\,f'(x_s)=0$ cannot hold
simultaneously. Hence acting by an automorphism (\ref{type-v}) with $q=z$ and a
general~$t$ does change the coordinate $v_s=0$, while keeping
nonzero those $v_i$ that were already nonzero. Now the claim
follows.

{\sl Step} 2. Suppose further that $v_i\neq 0$ for all
$i=1,\ldots,m$. Then acting with $G_v$ we can send our $m$-tuple to
another one where all the $u_i$, $i=1,\ldots,m$, are nonzero and
pairwise distinct. Indeed, let
$$
F(Q_i,q,t)=\frac{f'(x_i)}{1!}\frac{q(v_i)}{v_i}t+\ldots+
\frac{f^{(d)}(x_i)}{d!}\frac{q(v_i)^d}{v_i}t^d\in\kk[t]\,.
$$
We have $(x_i,v_i)\neq (x_j,v_j)$ for all $i\neq j$ because
$(x_i,u_i,v_i)\neq (x_j,u_j,v_j)$ where $u_i=f(x_i)/v_i$ and
$u_j=f(x_j)/v_j$. If $v_i=v_j$ then $f^{(d-1)}(x_i)\ne
f^{(d-1)}(x_j)$ since the linear form $f^{(d-1)}(x)$ is nonzero.
Thus for a suitable $q\in\kk[z]$ such that $q(v_i)\neq 0$ for all $i$,
the polynomials $F(Q_i,q,t)$ and $F(Q_j,q,t)$ are
different for every $i\ne j$. Applying the automorphism $H_v(q)$ from~(\ref{type-u})
with a general $t$, we obtain the result.

{\sl Step} 3. We assume now that all the coordinates $u_j$ are
nonzero and distinct. Let us show that it is possible, acting by $G_u$,
to map the initial tuple of points to a tuple having standard values of coordinates $v^{(0)}_s$, $s=1,\ldots,m$.
To this end, we construct an automorphism that preserves all the
points but $Q_i$ and sends $Q_i$ to a new point $Q'_i$ with
$v'_i=v^{(0)}_i$. Namely, fix a polynomial $q(z)$ with $q(0)=0$,
$q(u_i)\ne 0$ and $q(u_j)=0$ for all $j\ne i$. Our assumption on
$f(x)$ guarantees that the equation
$f(x)=u_i(v_i^{(0)}-v_i)+f(x_i)$ has a root $x=a_i$, where
$a_i\in\kk$. Applying $H_u(q)$ in ~(\ref{type-v}) with
$t=(a_i-x_i)/q(u_i)$, we obtain the required.

{\sl Step} 4. Suppose finally that $v_i=v_i^{(0)}$ for all $i$. It
suffices to reach the values $x_i=x_i^{(0)}$ for all $i$ acting by
an automorphism from $G_v$. Indeed, then also
$u_i=f(x_i)/v_i=f(x_i^{(0)})/v_i^{(0)}=u_i^{(0)}$. This can be
done by applying $H_u(q)$ as in~(\ref{type-u}) with $t=1$ and a
polynomial $q$ satisfying $q(0)=0$ and $q(v_i^{(0)})=x_i^{(0)} -
x_i$ for all~$i$. Now the proof is completed.
\end{proof}

\subsection{Infinite transitivity in higher dimensions}
It is enough to prove Theorems~\ref{main}, \ref{main10}
for $l=1$. Before passing to the proofs, we establish in
Lemmas~\ref{irr}-\ref{lemma1.2} some necessary elementary
facts concerning suspensions. In Lemmas~\ref{irr}-\ref{l15} we
consider an arbitrary field $\kk$ of characteristic zero.

\begin{lemma}\label{irr}
If $X^{(0)}$ is irreducible, then so is the suspension
$X^{(1)}=\SP(X^{(0)},f)$.
\end{lemma}
\begin{proof}
Suppose to the contrary that there exist nonzero functions
$$
F_1,F_2\in \kk[X^{(1)}]=\kk[X^{(0)}][u,v]/(uv-f) \qquad\text{such that}\quad
F_1F_2=0.
$$
We may assume that we chose $F_1$ and $F_2$ with minimal
$\deg_{u,v}(F_1)+\deg_{u,v}(F_2)$, and that no monomial in
$F_i$ contains the product $uv$, since otherwise we could replace
this product by $f$ according to Definition~\ref{susp}. If~$u$
occurs in both $F_1$ and $F_2$, then $\deg_u(F_1F_2)> 0$ since the
leading term in~$u$ cannot cancel. Hence, up to twisting~$u$ and~$v$,
we may assume that $F_1$ does not contain~$v$, and $F_2$ does
not contain~$u$. Let us write
$$
F_1=\sum_{i=0}^{k}a_iu^i, \qquad
F_2=\sum_{j=0}^{l}b_jv^j\,,
$$
where $a_i,b_j\in\kk[X^{(0)}]$, and $k+l$ is minimal possible.

If $k=l=0$, then $F_1$,
$F_2\in\kk[X^{(0)}]$ are zero divisors, which contradicts the
irreducibility of~$X^{(0)}$. So $k+l>0$.

If $a_0=b_0=0$, then we can decrease the degree $k+l$ by dividing out~$u$ and~$v$.
This contradicts the minimality assumption.
So we may assume that $a_0\neq 0$. Then the product $F_1F_2$
contains a nonzero term $a_0b_lv^l$, which gives again a
contradiction. The lemma is proved.
\end{proof}

\begin{lemma}\label{reg}
Let $\pi:X^{(1)}\to
X^{(0)}$ be the restriction of the projection $X^{(0)}\times\AA^2\to X^{(0)}$ to $X^{(1)}$.
Then $\pi (X^{(1)}_{\reg}) = X_{\reg}^{(0)}$.
\end{lemma}
\begin{proof}
Let $f_1, f_2,\ldots, f_m\in \kk[x_1,x_2,\ldots,x_s]$ be the functions generating the
ideal of $X^{(0)}\subseteq \AA^s$. A point $P\in X^{(0)}$ is
regular if and only if the rank of the Jacobian matrix
$$
D_0 = \qmatrix[1.3]{\frac{\partial f_1}{\partial x_1}
&\frac{\partial f_1}{\partial x_2}& \ldots & \frac{\partial
f_1}{\partial x_s}\\
\frac{\partial f_2}{\partial x_1}
&\frac{\partial f_2}{\partial x_2}& \ldots & \frac{\partial f_2}{\partial x_s}\\
\dotfill&\dotfill& \dotfill&\dotfill\\
\frac{\partial f_m}{\partial x_1}
&\frac{\partial f_m}{\partial x_2}& \ldots & \frac{\partial
f_m}{\partial x_s} }
$$
attains its maximal value $s-\dim X^{(0)}$ at~$P$. The
corresponding matrix for $X^{(1)}$ is
$$
D_1 = \qmatrix[1.3]{\frac{\partial f_1}{\partial x_1}
&\frac{\partial f_1}{\partial x_2}& \ldots & \frac{\partial
f_1}{\partial x_s}&0&0\\
\frac{\partial f_2}{\partial x_1}&\frac{\partial
f_2}{\partial x_2}& \ldots & \frac{\partial f_2}{\partial x_s}&0&0\\
\dotfill &\dotfill&\dotfill&\dotfill&\dotfill&\dotfill \\
\frac{\partial f_m}{\partial x_1} &\frac{\partial
f_m}{\partial x_2}& \ldots & \frac{\partial f_m}{\partial x_s}&0&0\\
\frac{\partial f}{\partial x_1} &\frac{\partial f}{\partial x_2}&
\ldots & \frac{\partial f}{\partial x_s}&-v&-u }.
$$
Obviously, $\rk D_1\leqslant 1+\rk D_0$ at every point. Since
$\dim X^{(1)}=1+\dim X^{(0)}$, any regular point of $X^{(1)}$ is
mapped via~$\pi$ to a regular point of $X^{(0)}$. On the other
hand, let $M$ be a square submatrix of $D_0$ and $P\in
X^{(0)}_{\reg}$ be a point such that $M(P)$ is of rank $r=s-\dim X^{(0)}$
equal to its order. We extend~$M$ to a square submatrix~$M'$ of
order $r+1$ by adding the last line and one of the two extra
columns of $D_1$ in such a way that $\rk M'(P,u,v)=1+\rk M(P)=r+1$
for some $(u,v)\neq (0,0)$, where $(P,u,v)\in X^{(1)}$. Then
$(P,u,v)\in X^{(1)}_{\reg}$. Now the assertion follows.
\end{proof}

\begin{remark}\label{sym}
Let~$A$ be an affine algebra. Recall (see~\cite{KZ}) that an {\em affine modification}
of~$A$ with center $(I,v)$, where $I\subseteq A$
is an ideal and $v\in I$ is not a zero divisor, is the quotient
algebra $A[It]/(1-vt)$, where
$$
A[It]=A\oplus \bigoplus_{n=1}^\infty (It)^n\cong
A\oplus I\oplus I^2\oplus \ldots = {\rm Bl}_I(A)
$$
is the Rees algebra of the pair $(A,I)$, and $t$ is a formal symbol.

Geometrically, the variety $\Spec\,(A[It]/(1-vt))$ can be obtained
from $X=\Spec A$ as follows. First we perform a blowup of~$X$ with the
center~$I$, and then remove the proper transform of
the zero divisor $V(v)$ in $X$ from ${\rm Bl}_I(X)$, which results again in an affine
variety. (We note that this proper transform meets the exceptional
divisor~$E$, since $v\in I$; see~\cite[\S 1]{KZ} for more
details).

According to~\cite[Example 1.4 and \S 5]{KZ}, the suspension
$X^{(1)}=\SP(X^{(0)},f)$ can be viewed as an affine modification
of $X^{(0)}\times \AA^1$ (where $\AA^1=\Spec \kk[v]$) with the center
$(I_1=(v,f),v)$ along the divisor $v=0$. Interchanging~$v$ and~$u$,
we may also regard the variety $X^{(1)}$ as an affine
modification of the product $X^{(0)}\times\AA^1$, where this time
$\AA^1=\Spec \kk[u]$, with the center $(I_2=(u,f),u)$ along the
divisor $u=0$. The exceptional divisors of these two modifications ($v=0$ and $u=0$,
respectively) are both isomorphic to
$X^{(0)}\times \AA^1$ but different as subvarieties of $X^{(1)}$.
For every $c\in\kk$ we denote by $U_c=\{u=c\}$ and $V_c=\{v=c\}$
the level hypersurfaces in $X^{(1)}$, which are widely used in the sequel.
\end{remark}

In~\cite[\S 2]{KZ} there was developed a method which allows to extend
an LND~$\partial$ to the affine modification provided that~$\partial$
stabilizes the center of the modification. In Lemma~\ref{lemma1.2}
(see below) we concretize this method in our particular case of
affine suspensions.

Given an LND~$\delta_0$ of an affine domain $A_{0}$ and a
polynomial $q\in\kk[z]$ with $q(0)=0$, we can define a new LND
$\delta'=\delta'(\delta_0,q)$ on $A'=A_0\otimes \kk[v]$, where~$v$
is a new variable, as follows. First we extend $\delta_0$ to~$A'$
by letting $\delta_0(v)=0$, and then we multiply $\delta_0$ by the
element $q(v)\in \ker \delta_0$. Suppose that $A_0$ is generated
by $x_1,x_2,\dots,x_s$. Then~$\delta'$ is given in coordinates by
\begin{equation}\label{a1}
\delta'(x_i)=q(v)\delta_0(x_i),\quad i=1,2,\ldots,s,\qquad
\delta'(v)=0.
\end{equation}

Let now~$u$ be yet another variable and $f\in A_0$ be nonzero.
Consider the structure algebra $A_1$ of the suspension over $A_0$:
$$
A_1=(A_0\otimes \kk[u,v])/(uv-f).
$$

\begin{lemma}\label{lemma1.2}
In the notation as above, every LND $\delta'\in\Der A'$ can be
transformed into an LND $\delta_1=\delta_1(\delta_0,
q)\in \Der A_1$ by letting
\begin{equation}\label{drob}
\delta_1(x_i)=\delta'(x_i),\quad i=1,2,\ldots,s\,,\qquad
\delta_1(u)= \frac{q(v)}{v}\,\delta_0(f) \,,\qquad
\delta_1(v)=\delta'(v)=0\,.
\end{equation}
\end{lemma}

\begin{proof}
First, let us check that these formulae extend~$\delta'$ to
$\delta_1=\delta_1(\delta_0, q)\in\Der (A_0\otimes \kk[u,v])$,
where~$\delta_1$ preserves
the ideal $(uv-f)$. Indeed, since $\frac{q(v)}{v}\in\kk[v]$ by our
choice of~$q$, the derivation $\delta_1$ is well defined on the
generators of $A_0\otimes \kk[u,v]$. It is easily seen that
$\delta_1$ is still locally nilpotent. The straightforward
calculation shows that $\delta_1(uv-f)=0$. Hence $\delta_1$
descends to an LND of the quotient algebra~$A_1$, denoted by the
same symbol $\delta_1$.
\end{proof}

\begin{definition}\label{gg1}
We let $G_v$ denote the subgroup of the special automorphism group
$\SAut\,(X^{(1)})$ generated by all one-parameter unipotent
subgroups
$$
H_v(\delta_0, q)=\exp(t\delta_1),\qquad\mbox{where}\quad
t\in\kk_+,\quad \delta_1=\delta_1(\delta_0,q)\,,
$$
with $\delta_0$ and $q(t)$ as above\footnote{Notice that for
$X^{(0)}=\AA^1=\Spec\kk[z]$ and $\delta_0=d/dz$ we have
$H_v(\delta_0, q)=H_v(q)$ from (\ref{type-u}).}.
Interchanging the
roles of~$v$ and~$u$, we obtain the second subgroup $G_u\subseteq
\SAut(X^{(1)})$. Here $u\in \kk[X^{(1)}]^{G_u}$ and $v\in
\kk[X^{(1)}]^{G_v}$. We will show that the subgroup $G\subseteq
\SAut\,(X^{(1)})$ generated by $G_u$ and $G_v$ acts infinitely
transitively in $X^{(1)}_{\reg}$.
\end{definition}

Given $k$ distinct constants $c_1,\ldots,c_k\in\kk$, we denote by
$\Stab^v_{c_1\ldots c_k}$ the subgroup of the group $G_v$
fixing pointwise all the hypersurfaces $V_{c_s}\subseteq X^{(1)}$,
$s=1,\ldots, k$.

\begin{lemma}\label{l15}
Suppose that the group $\SAut(X^{(0)})$ acts
$m$-transitively on $X^{(0)}_{\reg}$. Then for any distinct
constants $c_0,c_1,\ldots,c_k\in {\mathbb G}_{\rm m}(\kk)$ the group
$\Stab^v_{c_1\ldots c_k}$ acts $m$-transitively on $V_{c_0}\cap
X^{(1)}_{\reg}$.
\end{lemma}

\begin{proof}
Consider two collections of $m$ distinct points $P'_1,\ldots,P'_m$
and $Q'_1,\ldots,Q'_m$ in $V_{c_0}\cap X^{(1)}_{\reg}$. Denote by
$P_1,\ldots,P_m$ and $Q_1,\ldots,Q_m$ their
$\pi$-projections to $X^{(0)}$. Notice that the hypersurface
$V_{c_0}\subseteq X^{(1)}$ is mapped via $\pi$ isomorphically onto
$X^{(0)}$, while  by Lemma~\ref{reg} we have $\pi(V_{c_0}\cap
X^{(1)}_{\reg})\subseteq X^{(0)}_{\reg}$. A point $P'\in
V_{c_0}$ can be written as $P'=(P,u,c_0)$, where $P=\pi(P')\in
X^{(0)}$ and $u=u(P')=f(P)/c_0$. Consequently, every smooth point on $X^{(0)}$ has a
preimage on $X^{(1)}_{\reg}$. So these formulae give
an isomorphism $X^{(0)}\stackrel{\cong}{\longrightarrow} V_{c_0}$
which sends $P$ to~$P'$.

Since by our assumptions the group $\SAut(X^{(0)})$ acts
$m$-transitively on $X^{(0)}_{\reg}$, there exists an automorphism
$\psi_0\in \SAut(X^{(0)})$ which sends the ordered collection
$(P_1,\ldots,P_m)$ to $(Q_1,\ldots,Q_m)$. It can be written as a
product
$$
\psi_0=\prod_{i=1}^k\exp(\delta_0^{(i)})
$$
for some LNDs
$\delta_0^{(1)},\ldots,\delta_0^{(k)}\in\Der \kk[X^{(0)}]$.

Letting $q=\alpha z(z-c_1)\ldots(z-c_k)$, where
$\alpha\in {\mathbb G}_{\rm m}(\kk)$ is such that $q(c_0)=1$, by Lemma~\ref{lemma1.2}
we can lift all LNDs $\delta_0^{(i)}$ to the LNDs
$$
\delta_1^{(i)}=\delta_1^{(i)}(\delta_0^{(i)},q)
\in\Der\kk[X^{(1)}], \qquad i=1,\ldots,k\,.
$$
Respectively, $\psi_0$ can be lifted to an automorphism
$$
\psi_1=\prod_{i=1}^k\exp(\delta_1^{(i)})\in G_v\subseteq \SAut (X^{(1)})\,.
$$
By virtue of~$(\ref{a1})$ it is easily seen that the actions on
$X^{(1)}$ of the corresponding one-parameter unipotent subgroups
$H_v(\delta_0^{(i)},q)$ restrict to the original actions on
$V_{c_0}\cong X^{(0)}$. Hence the automorphism
$\psi_1\vert_{V_{c_0}}=\psi_0$ sends $(P'_1,\ldots,P'_m)$ to
$(Q'_1,\ldots,Q'_m)$. Due to our choice of~$q(z)$, this
automorphism fixes all the other hypersurfaces $V_{c_s}$
pointwise. The lemma is proved.
\end{proof}

\begin{lemma}\label{l16}
Let $\kk$ be an algebraically closed
field of characteristic  zero. Suppose as before that the group
$\SAut(X^{(0)})$ acts $m$-transitively on $X^{(0)}_{\reg}$. Then
for any set of distinct points $Q'_1,\ldots,Q'_m\in X^{(1)}_{\reg}$
there exists an automorphism $\varphi\in\SAut (X^{(1)})$ such
that $\varphi(Q'_i)\not\in U_0\cup V_0$ for all $i=1,2,\ldots,m$.
\end{lemma}

\begin{proof}
We say that the point $Q'_i=(Q_i, u_i,v_i)\in X^{(1)}$ is {\it
hyperbolic} if $u_iv_i\ne 0$, i.e., if $Q'_i\not\in U_0\cup V_0$. We
have to show that the original collection can be moved by means of
an automorphism from $\SAut (X^{(1)})$ to a new collection so that all the points become hyperbolic.
Suppose that $Q'_1,\ldots,Q'_l$ are already hyperbolic while
$Q'_{l+1}$ is not, where $l\ge 0$. By recursion, it is sufficient
to move $Q'_{l+1}$ off $U_0\cup V_0$ while leaving the points
$Q'_1,\ldots,Q'_l$ hyperbolic.
Consider the following two cases:

{\sl Case} 1: $u_{l+1}=0$, $v_{l+1}\ne 0$.

{\sl Case} 2: $u_{l+1}=v_{l+1}=0$.

We claim that there exists an automorphism
$\varphi\in\SAut(X^{(1)})$ leaving $Q'_1,\ldots,Q'_l$ hyperbolic
such that  in Case~1 the point $\varphi(Q'_{l+1})$ is hyperbolic
as well, and  in Case~2 this point  satisfies the assumptions of
Case~1.

In Case~1 we split $Q'_1,\ldots,Q'_{l+1}$ into several disjoint
pieces $M_0,\ldots,M_k$ according to different values of~$v$:
$Q'_i\in M_j\iff v_i=c_j$, where $c_j\neq 0$.
Assume that $M_0=\{Q'_{i_1}, \ldots, Q'_{i_r}, Q'_{l+1}\}$,
where $i_k\le l$ for all $k=1,\ldots,r$. We can choose an extra
point $Q''_{l+1}\in (V_{c_0}\cap X^{(1)}_{\reg} ) \setminus U_0$.
Indeed, since $c_0=v_{l+1}\neq 0$, we have $V_{c_0}\cong X^{(0)}$.
Under the assumptions of Theorem~\ref{main} $\dim X^{(0)}\ge 2$,
hence also $\dim ((V_{c_0}\cap X^{(1)}_{\reg} ) \setminus U_0)=\dim
X^{(0)}\ge 2$.

By Lemma~\ref{l15} the subgroup
$\Stab^v_{c_1,\ldots,c_{k}}\subseteq G_v$ acts
$(r+1)$-transitively on $V_{c_0}\cap X^{(1)}_{\reg} $. Therefore
we can send the $(r+1)$-tuple $(Q'_{i_1},\ldots, Q'_{i_r},
Q'_{l+1})$ to $(Q'_{i_1},\ldots, Q'_{i_r}, Q''_{l+1})$ fixing the
remaining points of $M_1\cup \ldots \cup M_{k}$. This confirms our
claim in Case~1.

In Case~2 the point $Q'_{l+1}=(Q_{l+1},0,0)$ belongs to $X^{(1)}_{\reg}$.
From Lemma~\ref{reg} and its proof it follows that
$Q_{l+1}=\pi(Q'_{l+1})\in X^{(0)}_{\reg}$ and $df(Q_{l+1})\neq 0$
in the cotangent space $T^*_{Q_{l+1}} X^{(0)}$. The variety
$X^{(0)}$ being flexible, there exists an LND $\partial_0\in\Der
\kk[X^{(0)}]$ such that $\partial_0(f)(Q_{l+1})\neq 0$. Letting
$q(v)=v(v-v_1)(v-v_2)\ldots(v-v_l)\in\kk[v]$ and choosing a set of
generators $x_1,\ldots,x_s$ of the algebra $\kk[X^{(0)}]$,
similarly as in~(\ref{drob}) we extend $\partial_0$ to
$\partial_1\in \Der\, \kk[X^{(1)}]$ via
\begin{equation}\label{drob1}
\partial_1 (x_i)  = q(v)\partial_0 (x_i),\quad i=1,2,\ldots, s\,,\qquad
\partial_1 (u) = \frac{q(v)}v\, \partial_0 (f)\,,\qquad
\partial_1 (v) = 0\,.
\end{equation}
Due to our choice, $\partial_1 (u)(Q'_{l+1})\neq 0$. Hence the
action of the associate one-parameter unipotent subgroup
$H_v(\partial_0,q)=\exp (t\partial_1)$ pushes the point $Q'_{l+1}$
out of~$U_0$. So the orbit $H_v(\partial_0,q).Q'_{l+1}$ meets the
hypersurface $U_0\subseteq X^{(1)}$ in finitely many points.
Similarly, for every $j=1,2,\ldots,l$ the orbit
$H_v(\partial_0,q).Q_j\not\subseteq U_0$ meets $U_0$ in finitely
many points. Letting now $\varphi=\exp(t_0
\partial_1)\in H_v(\partial_0,q)\subseteq G_v$, we conclude that for
a general value of $t_0\in\kk$ the images $\varphi(Q'_j)$ lie
outside~$U_0$ for all $j=1,2,\ldots, l+1$. Since the group
$H_v(\partial_0,q)$ preserves the coordinate~$v$, the points
$\varphi(Q'_1),\ldots,\varphi(Q'_l)$ are still hyperbolic.
Interchanging\footnote{We have not done
this earlier in order to keep our previous notation.} the role of~$u$ and~$v$  we see that the assumptions of Case~1 are fulfilled for the new
collection $\varphi(Q'_1),\ldots,\varphi(Q'_l),\varphi(Q'_{l+1})$,
as required.
\end{proof}

{\sc Proof of infinite transitivity in
Theorem~\ref{main}.}
If $X^{(0)}=\AA^1$ then the assertion follows
from Theorem~\ref{aa}. Let now $\dim X^{(0)}\ge 2$. To show that
the action of the group $\SAut(X^{(1)})$ on $X^{(1)}_{\reg}$ is
$m$-transitive for every $m\in\NN$, we fix a standard collection
of~$m$ distinct points $P'_1,\ldots,P'_m\in U_1\cap
X^{(1)}_{\reg}$. It suffices to show that any other $m$-tuple of
distinct points $Q'_1,\ldots,Q'_m\in X^{(1)}_{\reg}$ can be mapped to $P'_1,\ldots,P'_m$ by means of an automorphism
$\psi\in \SAut(X^{(1)})$. In view of Lemma~\ref{l16} we may
suppose that $Q'_i\not\in U_0\cup V_0$ for all $i=1,\ldots,m$.
Similarly as in the proof of Lemma~\ref{l16} we split the
collection $Q'_1,\ldots,Q'_m$ into disjoint pieces $M_1,\ldots,
M_k$ according to the values of coordinate~$v$.

By our assumption the variety $X^{(0)}$ is flexible. It follows
that the only units in $\kk[X^{(0)}]$ are constants. Consequently,
since~$f$ is non-constant we have $f(X^{(0)})=\kk$. In particular,
$U_c\cap V_d\neq\varnothing$ for any $c,d\in\kk$. Since $\dim
X^{(1)}=1+\dim X^{(0)}\ge 3$, the intersection $U_c\cap V_d$ has
positive dimension, hence is infinite.

Therefore, acting with the subgroups $\Stab^v_{c_1\ldots\check c_i
\ldots c_l}\subseteq G_v$, by Lemma~\ref{l15} we can send~$M_i$ to
$U_1 \cap V_{c_i}\cap X^{(1)}_{\reg}$ fixing all the other points from
$\bigcup_{j\ne i} M_j$.
So we may assume that
$Q'_1,\ldots,Q'_m\in U_1\cap X_{\reg}$.  Applying Lemma~\ref{l15}
again with~$u$ and~$v$ interchanged, $k=0$,
and $c_0=1$, i.e., acting with the subgroup $G_u$,
we can send the resulting collection to the standard one
$P'_1,\ldots, P'_m$. Now the proof is completed.


\subsection{Suspensions over real varieties}
Here we prove
Theorem~\ref{main10}. We need the following elementary lemma.

\begin{lemma}\label{level}
Let $Y$ be a smooth connected real
manifold of dimension at least~2. Then for any continuous
function $f:Y \to\RR$ the level set $f^{-1}(c)$ is infinite for
each $c\in\Int f(Y)$.
\end{lemma}

\begin{proof}
Under our assumptions $\Int f(Y)\subseteq\RR$
is an open interval. Choose two points $y_1,y_2\in Y$ such that
$f(y_1)=c_1<c$ and $f(y_2)=c_2>c$. We can join them by a
smooth path~$l$ in~$Y$. There exists a tubular neighborhood~$U$ of~$l$
diffeomorphic to a cylinder $\Delta\times I$, where $I=[0,1]$ and
$\Delta$ is a ball of dimension $\dim \Delta=\dim Y-1\ge 1$. So
there exists a continuous family of paths joining~$y_1$ and~$y_2$
within~$U$ such that any two of them meet only at their ends~$y_1$
and~$y_2$. Since the hypersurface $f^{-1}(c)$ separates~$Y$, each
of these paths crosses it. In particular, $f^{-1}(c)$ is infinite.
\end{proof}

The proof of Theorem~\ref{main10} differs just slightly from that
of Theorem~\ref{main}. Hence it is enough to indicate the
necessary changes.

\medskip

{\sc Sketch of the proof of Theorem~\ref{main10}.}
The assumption that the field~$\kk$ is algebraically closed was
actually  used in the proof of Theorem~\ref{main} only on two
occasions. Namely, in the proofs of Lemma~\ref{l16} and of the
infinite transitivity in Theorem~\ref{main} we exploited the fact
that under our assumptions the level sets $(V_{c_k}\cap
X^{(1)}_{\reg} )\setminus U_0$ and $U_1 \cap V_{c_i}\cap
X^{(1)}_{\reg} $ are of positive dimension, hence are infinite.
For $\kk=\overline{\kk}$ the latter follows from the Krull theorem
and the dimension count. In the case where $\kk=\RR$, we can
deduce the same conclusion using Lemma~\ref{level}. Indeed, in the
notation as above, for every $c_i\ne 0$ the restrictions
$$
\pi\colon V_{c_i}\cap X^{(1)}_{\reg}=(V_{c_i})_{\reg} \to
X_{\reg} ^{(0)},\qquad \pi\colon U_1 \cap V_{c_i}\cap
X^{(1)}_{\reg} \to f^{-1}(c_i)\cap X_{\reg} ^{(0)}
$$
are isomorphisms. Under the assumptions of Theorem~\ref{main10} the
smooth real manifold $X_{\reg} ^{(0)}$ has dimension $\ge 2$ and is
connected. Since $f(X_{\reg} ^{(0)})=\RR$, by Lemma~\ref{level}
the level set $f^{-1}(c_i)\cap X_{\reg} ^{(0)}$ is infinite. Since
$c_k\ne 0$, the set $(V_{c_k}\cap X^{(1)}_{\reg} )\setminus
U_0\supseteq U_1 \cap V_{c_k}\cap X^{(1)}_{\reg}$ is infinite too.

A posteriori, the manifold $X^{(1)}_{\reg}$ is also connected.
Hence by recursion we can apply the argument to the iterated
suspensions $X^{(i)}$ over $X^{(0)}$, $i=1,\ldots,l$.

\subsection{Flexibility}
To complete the proofs of Theorems~\ref{aa}--\ref{main10}, it
remains to establish the flexibility of $X^{(1)}$.

\begin{lemma}\label{lastl}
Under the assumptions
of any one of Theorems~\ref{aa}--\ref{main10} the variety $X^{(1)}$
is flexible.
\end{lemma}

\begin{proof}
We already know that the group $\SAut(X^{(1)})$
acts transitively on $X_{\reg}^{(1)}$. Hence, similarly as in the
proof of Lemma~\ref{flex}, it suffices to find just one flexible
point $P'=(P,u,v)\in X_{\reg}^{(1)}$.

The function $f\in\kk[X^{(0)}]$ being non-constant,  $df(P)\neq 0$
at some point $P\in X_{\reg}^{(0)}$ with $f(P)\neq 0$. Due to our
assumption $X^{(0)}$ is flexible. Hence there exist $n$ locally
nilpotent derivations
$\partial_0^{(1)},\ldots,\partial_0^{(n)}\in\Der\kk[X^{(0)}]$,
where $n=\dim X^{(0)}$, such that the corresponding vector fields
$\xi_1,\ldots, \xi_n$ span the tangent space $T_{P}X^{(0)}$, i.e.,
$$
\rk \left(%
\begin{array}{c}
 \xi_1(P) \\
 \ldots \\
 \xi_n(P) \\
\end{array}%
\right)=n\,.
$$
It follows that $\partial_0^{(i)} (f)(P)\ne 0$ for some index
$i\in\{1,\ldots,n\}$.

Let now $P'=(P,u_0,v_0)\in X_{\reg}^{(1)}$ be a point such that
$\pi(P')=P$. Since $u_0v_0=f(P)\neq 0$, the point $P'$ is
hyperbolic. Letting $q(v)=v$ in Lemma~\ref{lemma1.2} we obtain $n$
LNDs
$$
\partial_1^{(1)},\ldots,\partial_1^{(n)}\in\Der\kk[X^{(1)}],\quad\mbox{
where}\quad \partial_1^{(j)}=\partial_1^{(j)}(\partial_0^{(j)},v)\,.
$$
Interchanging~$u$ and~$v$ and letting $j=i$, we get yet another LND
$$
\partial_2^{(i)}=\partial_2^{(i)}(\partial_0^{(i)},u) \in\Der\kk[X^{(1)}]\,.
$$
Let us show that the corresponding $n+1$ vector fields span the
tangent space $T_{P'}X^{(1)}$ at~$P'$, as required. We can view
$\partial_1^{(1)},\ldots,\partial_1^{(n)},\partial_2^{(i)}$ as
LNDs in $\Der\kk[X^{(0)}][u,v]$ preserving the ideal $(uv-f)$, which means
that the corresponding vector fields are tangent to the
hypersurface
$$
X^{(1)}=\{uv-f(P)=0\}\subseteq X^{(0)}\times\AA^2\,.
$$
The values  of these vector fields at
$P'\in X^{(1)}_{\reg}$ yield an $(n+1)\times (n+2)$-matrix
$$
E=\left(%
\begin{array}{ccc}
 v_0\xi_1(P) & \partial_0^{(1)}(f)(P) & 0 \\
 \dotfill & \dotfill & \dotfill \\
 v_0\xi_n(P) & \partial_0^{(n)}(f)(P) & 0 \\
  u_0\xi_i(P) & 0 & \partial_0^{(i)}(f)(P) \\
\end{array}%
\right)\,.
$$
The first~$n$ rows of~$E$ are linearly independent, and the last
one is independent from the preceding ones since
$\partial_0^{(i)}(f)(P)\neq 0$. Therefore, $\rk (E)=n+1=\dim
X^{(1)}$. So our locally nilpotent vector fields indeed span the
tangent space $T_{P'}X^{(1)}$ at $P'$, as claimed.
\end{proof}

Now the proofs of Theorems~\ref{aa}-\ref{main10} are completed.




\begin{thebibliography}{}
%
\bibitem{Ak} D.N.\ Ahiezer (Akhiezer): Dense orbits with two endpoints.
Izv.\ Akad.\ Nauk SSSR Ser.\ Mat.\ 41 (1977), 308--324 (Russian). English
transl.: Math.\ USSR-Izv.\ 11 (1977), 293--307 (1978).
%
\bibitem{AL} E.\ Anders\'en and L.\ Lempert: On the Group of Holomorphic
Automorphisms of $\CC^n$. Invent.\ Math.\ 110 (1992), 371--388.
%
\bibitem{AFKKZ}
I.\ Arzhantsev, H.\ Flenner, S.\ Kaliman, F.\ Kutzschebauch, and M.\ Zaidenberg:
Flexible varieties and automorphism groups.
arXiv:1011.5375 (2010), 41 p, Duke Math. J. (to appear).
%
\bibitem{BH} I.\ Biswas and J.\ Huisman:
Rational real algebraic models of topological surfaces.
Doc.\ Math.\ 12 (2007), 549--567.
%
\bibitem{BM} J.\ Blanc and F.\ Mangolte:
Geometrically rational real conic bundles and very transitive
actions. Compos. Math. 147 (2011), 161--187.
%
\bibitem{Bo} A.\ Borel: Les bouts des espaces homog\`enes de groupes de
Lie.
Ann.\ Math.\ (2) 58 (1953), 443--457.
%
\bibitem{CLS} D.A.\ Cox, J.B.\ Little, and H.\ Schenck:
Toric Varieties. Graduate Studies in Math., Vol.~124, Amer. Math. Soc.,
Providence, Rhode Island, 2011.
%
\bibitem{Da} D.\ Daigle:
On locally nilpotent derivations of $k[X_1 , X_2 , Y ]/(\varphi(Y
) - X_1 X_2 ) $.  J.\ Pure Appl.\ Algebra  181 (2003), 181--208.
%
\bibitem{Dem}
M.\ Demazure: Sous-groupes alg\'{e}briques de rang maximum du
groupe de Cremona. Ann. Sci. \'Ecole Norm. Sup. (4) {3} (1970),
507--588.
%
\bibitem{FKZ3}  H.\ Flenner, S.\ Kaliman, and M.\ Zaidenberg:
Uniqueness of $\CC^*$- and $\CC_+$-actions on Gizatullin surfaces.
Transform.\ Groups 1 (2008), 305--354.
%
\bibitem{FZ2}
H.\ Flenner and M.\ Zaidenberg: Locally nilpotent derivations on
affine surfaces with a $\CC^*$-action. Osaka J.\ Math. {42}
(2005), 931-074.
%
\bibitem{Fo} F.\ Forstneric: Interpolation by holomorphic automorphisms
and
embeddings in $\CC^n$. J.\ Geom.\ Anal.\ 9 (1999), 93--117.
%
\bibitem{Fre}
G.\ Freudenburg: Algebraic Theory of Locally Nilpotent
Derivations. Invariant Theory and Algebraic Transformation groups, vol. VII, Encyclopaedia of Mathematical Sciences, Vol.~{136}, Springer-Verlag, Berlin 2006.
%
\bibitem{Fu}
W.\ Fulton: Introduction to toric varieties. Annals of Math.
Studies {131}, Princeton University Press, Princeton, NJ, 1993.
%
\bibitem{Gi}
M.H.\ Gizatullin: Affine surfaces which are quasihomogeneous with
respect to an algebraic group.
Izv.\ Akad.\ Nauk SSSR Ser.\ Mat.\  35:4 (1971), 738 -- 753 (Russian);
English transl.: Mat. USSR Izv. {5}:4 (1971), 754-769.
%
\bibitem{HO} A.\ Huckleberry and A.\ Oeljeklaus: A characterization of
Complex
Homogeneous Cones. Math.\ Zeit.\ 170 (1980), 181--194.
%
\bibitem{HM1}
J.\ Huisman and F.\ Mangolte: The group of automorphisms of a real
rational surface is $n$-transitive. Bull.\ Lond.\ Math.\ Soc.\ 41
(2009), 563--568.
%
\bibitem{HM2}
J.\ Huisman and F.\ Mangolte: Automorphisms of real rational
surfaces and weighted blow-up singularities. Manuscripta \ Math.\ 132 (2010), 1--17.
%
\bibitem{KK1} S.\ Kaliman and F.\ Kutzschebauch:
Density property for hypersurfaces $UV=P(\overline X)$. Math.\
Zeit.\ 258 (2008), 115--131.
%
\bibitem{KK2} S.\ Kaliman and F.\ Kutzschebauch: On the present state of
the
Andersen-Lempert theory. CRM Proceedings and Lecture Notes, Vol.~54,
Amer. Math. Soc. (2011), 85--122.
%
\bibitem{KZ}
S.\ Kaliman and M.\ Zaidenberg: Affine modifications and affine
hypersurfaces with a very transitive automorphism group.
Transform.\ Groups 4 (1999), 53--95.
%
\bibitem{KPZ}
T.\ Kishimoto, Yu.\ Prokhorov, and M.\ Zaidenberg: Group actions
on affine cones. CRM Proceedings and Lecture Notes, Vol.~54,
Amer. Math. Soc. (2011), 123--164.
%
\bibitem{Kr} L.\ Kramer: Two-transitive Lie groups. J.\ Reine Angew.\
Math.\ 563 (2003), 83--113.
%
\bibitem{KM}
K.\ Kuyumzhiyan and F.\ Mangolte:
Infinitely transitive actions on real affine suspensions.
J. Pure Appl. Algebra, {216}:10 (2012), 2106-2112.
%
\bibitem{LR} V.\ Lakshmibai, K.N.\  Raghavan: Standard Monomial Theory. Invariant Theory and Algebraic Transformation groups, vol. VIII, Encyclopaedia of Mathematical Sciences, Vol.~{137}, Springer-Verlag, Berlin 2008
%
\bibitem{Le} F.\ Lescure: \'Elargissement du groupe d'automorphismes
pour des
vari\'et\'es quasi-homogenes. Math.\ Ann.\ 261 (1982), 455--462.
%
\bibitem{Li}
A.\ Liendo: Affine $T$-varieties of complexity one and locally
nilpotent derivations. Transform.\ Groups\ 15 (2010), no.~2, 389–425.
%
\bibitem{Lie2} A.\ Liendo:
${\mathbb G}_{\rm a}$-actions of fiber type on affine $T$-varieties. J.\ Algebra
324 (2010), 3653--3665.
%
\bibitem{ML1} L.\ Makar-Limanov: On groups of
automorphisms of a class of surfaces. Israel J.\ Math.\ 69 (1990),
250--256.
%
\bibitem{ML2} L.\ Makar-Limanov:
Locally nilpotent derivations on the surface $xy=p(z)$.
Proceedings of the Third International Algebra Conference (Tainan,
2002),  215--219. Kluwer Acad.\ Publ.\, Dordrecht, 2003.
%
\bibitem{Perep}
A.\ Perepechko: Flexibility of affine cones over del Pezzo surfaces of degree 4 and 5.
arXiv:1108.5841 (2011), 6~p., to appear in: Func. Analysis and Its Appl.
%
\bibitem{Po}
V.L.\ Popov: Classification of affine algebraic surfaces that are
quasihomogeneous with respect to an algebraic group.
Izv.\ Akad.\ Nauk SSSR Ser.\ Mat.\  37:5 (1973), 1038--1055 (Russian);
English transl.: Math.\ USSR
Izv.\ 7:5 (1973), 1039-1055 (1975).
%
\bibitem{Po1}
V.L.\ Popov: Picard groups of homogeneous spaces of linear
algebraic groups and one-dimensional homogeneous vector fiberings.
Izv.\ Akad.\ Nauk SSSR Ser.\ Mat.\ 38 (1974), 294--322 (Russian);
English transl.: Math.\ USSR-Izv.\ 8 (1974), 301--327 (1975).
%
\bibitem{Po3}
V.L.\ Popov: Generically multiple transitive algebraic group
actions. Algebraic groups and homogeneous spaces,  481--523, Tata
Inst.\ Fund.\ Res.\ Stud.\ Math., Mumbai, 2007.
%
\bibitem{Po2}
V.L.\ Popov: On the Makar-Limanov, Derksen invariants, and finite
automorphism groups of algebraic varieties. CRM Proceedings and Lecture Notes, Vol.~54, Amer. Math. Soc. (2011), 289--311.
%
\bibitem{PV}
V.L.\ Popov and E.B.\ Vinberg: A certain class of quasihomogeneous
affine varieties. Izv.\ Akad.\ Nauk SSSR Ser.\ Mat.\ 36 (1972),
749--764 (Russian); English transl.: Math.\ USSR Izv.\ 6 (1972),
743--758.
%
\bibitem{PV2}
V.L.\ Popov and E.B.\ Vinberg: Invariant Theory. In: Algebraic
geometry IV,  A.N.\ Parshin, I.R.\ Shafarevich (eds), Berlin,
Heidelberg, New York: Springer-Verlag, 1994.

%
\bibitem{Ro} J.-P.\ Rosay: Automorphisms of $\CC^n$, a survey
of Anders\'en-Lempert theory and applications. Complex geometric
analysis in Pohang (1997), 131--145, Contemp. Math., 222, Amer.\
Math.\ Soc., Providence, RI, 1999.
%
\bibitem{RR} J.-P.\ Rosay and W.\ Rudin:
Holomorphic maps from $\CC^n$ to $\CC^n$.
Trans.\ Amer.\ Math.\ Soc.\ 310 (1988), 47--86.
%
\bibitem{Ti} J.\ Tits:
Sur certaines classes d'espaces homog\`enes de groupes de Lie.
Acad.\ Roy.\ Belg.\ Cl.\ Sci.\ M\'em.\ Coll.\ in $8^\circ$, 29
(1955), 268 pp.
%
\bibitem{TV} \'A.\ T\'oth and D.\ Varolin: Holomorphic diffeomorphisms
of
semisimple homogeneous spaces. Compos.\ Math.\ 142 (2006),
1308--1326.
%
\bibitem{Va} D.\ Varolin: The density property for complex manifolds and
geometric structures. II. Internat.\ J.\ Math.\ 11 (2000),
837--847.
%
\end{thebibliography}
\end{document}